%% file: IWFPM.tex
\documentclass[review,onefignum,onetabnum,nohypdvips]{siamart171218}

\input{ex_shared}

\begin{document}

\maketitle
	
\begin{abstract}
    In this paper, we propose a new algorithm, the irrational-window-filter projection method (IWFPM), for quasiperiodic systems with concentrated spectral point distribution. 
    Based on the projection method (PM), IWFPM filters out dominant spectral points by defining an irrational window and uses a corresponding index-shift transform to make the FFT available. 
    The error analysis on the function approximation level is also given. 
    We apply IWFPM to 1D, 2D, and 3D quasiperiodic Schr\"odinger eigenproblems (QSEs) to demonstrate its accuracy and efficiency. 
    IWFPM exhibits a significant computational advantage over PM for both extended and localized quantum states. 
    More importantly, by using IWFPM, the existence of Anderson localization in 2D and 3D QSEs is numerically verified.
\end{abstract}



\section{Introduction}
\label{sec:int}

Quasiperiodic systems, as a natural extension of periodic structures, have been widely observed in physics and materials sciences, such as many-body problems, quasicrystals, incommensurate systems, polycrystalline materials, and quantum systems~\cite{poincare1890probleme,shechtman1984metallic,cao2021computing,sutton1995interfaces,hasan2010colloquium}.
Over these years, a growing realization has emerged that underlying irrational numbers of quasiperiodic systems impart various fascinating features~\cite{penrose1974role,meyer2000algebraic,baake2013aperiodic,bohr2018almost}. 
Particularly, in quantum systems, numerous intriguing physical phenomena have been discovered to be related to quasiperiodic structures, including quantum Hall effect, Anderson localization, topological insulators, photonic moir\'e lattices, and mobility edge~\cite{thouless1982quantized,hofstadter1976energy,hasan2010colloquium,wang2020localization,sarma1988mobility,wang2020one}.

Quasiperiodic systems pose significant challenges for numerical simulations, due to their space-filling order without decay or translation invariance. 
Recently years, several methods for solving quasiperiodic systems have been developed. 
The widely used periodic approximation method \cite{zhang2008efficient} employs periodic solutions to approximate quasiperiodic solutions,  
inevitably introducing rational approximation errors~\cite{jiang2014numerical,jiang2023on}. 
An accurate algorithm is the projection method (PM), which treats the quasiperiodic system as an irrational manifold of a high-dimensional periodic system~\cite{jiang2014numerical}. 
PM has spectral accuracy, and is efficient owing to its utilization of fast Fourier transform (FFT)~\cite{jiang2024numerical}. 
Further, the finite points recovery method is proposed for both high- and low-regularity quasiperiodic systems~\cite{jiang2024accurately}.

\textbf{Motivation}.
PM has shown outstanding advantages in accurately computing quasiperiodic systems, especially in incommensurate quantum systems~\cite{jiang2024numerical,li2021numerical,wang2022effective,jiang2023high,zhou2019plane}.
However, when using PM, quasiperiodic function after lifting dimension may exhibit distinct regularities along different directions, leading to a notable deterioration in convergence. 
This flaw becomes apparent when solving some quasiperiodic systems with singularity solution, like localized quantum state.

Focusing on Anderson localization, the phenomenon of wave diffusion being absent in a disordered or quasiperiodic medium, it is of great significance in regulating various physical properties in materials, including conductivity, optical properties, and magnetism~\cite{anderson1958absence}. 
Over the past few decades, incommensurate electrical structures, drawing attention for their capacity of achieving the continuous transition from extended state to localized state, have been experimentally studied through techniques like cold atom control and optical superlattice~\cite{roati2008anderson,deissler2010delocalization,stutzer2018photonic,wang2020localization}. 
Under the tight-binding limit, the Hamiltonian of incommensurate quantum system can be mapped onto the well-known almost Mathieu operator in 1D discrete case, which is a typical form of quasiperiodic Schr\"{o}dinger eigenproblems (QSEs)~\cite{simon2000schrodinger}. 
Since the 1980s, substantial progress has been made in the spectral theory of QSEs. 
Researchers have found that the spectral structure of QSEs can be decomposed into pure point, singular continuous, and absolutely continuous spectra, which correspond to localized, critical, and extended states of quantum systems, respectively~\cite{bellissard1982cantor, bourgain2002continuity,simon2000schrodinger,bourgain2007anderson,avila2009ten,avila2015global,avila2017sharp,ge2023multiplicative,shi2021absence}.
While significant theoretical works have been conducted on the one-dimensional cases of QSEs, addressing two- and higher-dimensional scenarios remains a challenging endeavor~\cite{bourgain2007anderson,shi2021absence}. 

When numerically solving arbitrary dimensional QSEs, especially in these cases where the wavefunction exhibits localization, PM might become inefficient due to the high computational cost. 
Based on PM, several heuristic works have been introduced. 
Ref. \cite{wang2022convergence} observed the phenomenon of concentrated distribution of spectral points in QSEs. 
Based on this observation, Refs. \cite{wang2022convergence,gao2023reduced} both capture the concentrated distribution of spectral points by a parallelogram index set, thus reducing the degrees of freedom of PM.
Unfortunately, the RPM does not improve the FFT computational efficiency of PM, as it adopts a zero-fill operation when 
performing FFT on such an irregular index set. 
How to improve PM to make the FFT available and applied to QSEs to find high-dimension Anderson localization is the main purpose of this paper. 

\textbf{Contribution.}
In this paper, we propose a new algorithm, named the \textit{irrational-window-filter projection method} (IWFPM), and apply it to arbitrary dimensional global quasiperiodic systems. 
Based on the PM and the phenomenon that the spectral points are concentrated along an irrational direction,  
IWFPM filters out dominant spectral points by defining an irrational window.
Moreover, a corresponding index-shift transform is designed to make the FFT available.  
The error analysis on the function approximation level is also given.
We apply IWFPM to 1D, 2D, and 3D QSEs to demonstrate its accuracy and efficiency. 
An efficient diagonal preconditioner is also designed for the discrete QSEs to significantly reduce condition number. 
Numerous experiments demonstrate that IWFPM has an absolute computational advantage over PM for both extended and localized quantum states. 
More importantly, by using IWFPM, the existence of Anderson localization in 2D and 3D QSEs is numerically verified.

\textbf{Organization.}
The article is structured as follows. In \cref{sec:pre}, we introduce the preliminaries about quasiperiodic functions and give a brief introduction of the PM. 
In \cref{sec:IWFPM}, we present the IWFPM and its implementation process. 
Moreover, we define a new norm in quasiperiodic function space to enable the convergence analysis of IWFPM. 
In \cref{sec:num}, we illustrate the effectiveness and superiority of IWFPM through its application to 1D, 2D, and 3D QSEs, and verify the existence of Anderson localization. 
Finally, in \cref{sec:concl}, we summarize this work and give an outlook on future work.

\section{Quasiperiodic functions and projection method (PM)}\label{sec:pre}
In this section, we present the preliminaries about quasiperiodic functions and offer a brief overview of PM. 
\begin{definition}\label{def:pro_mat}
    A matrix $\bm{P}\in\mathbb{R}^{d\times n}$ is called the projection matrix, if it belongs to the set
    $    \mathbb{P}^{d\times n} :=\{\bm{P}=(\bm{p}_1,\cdots,\bm{p}_n)\in\mathbb{R}^{d\times n}:
    \bm{p}_1,\cdots,\bm{p}_n~\mbox{are}~\mathbb{Q}\mbox{-linearly}~\mbox{independent},\\ {\rm rank}(\bm{P})=d\}.$
\end{definition}
\begin{definition} \label{df:quasi}
	A $d$-dimensional function $u(\bm{x})$ is quasiperiodic if there exists a continuous $n$-dimensional periodic function ${U}(\bm{y})$ and a projection matrix $\bm{P}\in \mathbb{P}^{d\times n}$, such that 
$    u(\bm{x})={U}(\bm{P}^T \bm{x})$, for all    $\bm{x}\in\mathbb{R}^d.$
\end{definition}

\begin{remark}
    The periodic function ${U}(\bm{y})$ is called the parent function of $u(\bm{x})$. 
    And we use the notation $\mathcal{Q}(\mathbb{R}^d)$ to represent the set of all $d$-dimensional quasiperiodic functions. Without loss of generality, we always assume that all parent functions are measurable on $n$-dimensional torus $\mathbb{T}^n:=(\mathbb{R}/2\pi\mathbb{Z})^n$.
\end{remark}

For $n$-dimensional periodic functions ${U}$ and ${V}$, their inner product is 
\begin{equation*}\label{eq:inner}
    \left<{U},{V}\right>:=\frac{1}{(2\pi)^n}\int_{[0,2\pi]^n}{U}(\bm{y})\overline{{V}}(\bm{y})\mathrm{d}\bm{y}.
\end{equation*}
We say ${U}\in\mathcal{L}^2(\mathbb{T}^n)$ if
$\|{U}\|_{\mathcal{L}^2}:=\left<{U},{U}\right>^{1/2}<+\infty$. 
Denote Fourier basis function
$\varphi_{\bm{k}}(\bm{y}):=\mathrm{e}^{i\bm{k\cdot y}}$
for index $\bm{k}\in\mathbb{Z}^n$, where $\bm{k\cdot y}=\sum_{j=1}^nk_jy_j$. It is obvious that for any $\bm{k},~\bm{k}'\in\mathbb{Z}^d$, the Fourier basis functions $\varphi_{\bm{k}}$ and $\varphi_{\bm{k}'}$ are orthogonal, i.e.
\begin{equation*}
    \left<\varphi_{\bm{k}},\varphi_{\bm{k}'}\right>:=\delta_{\bm{k}\bm{k}'}=
    \left\{\begin{array}{ll}
    1, & \bm{k}=\bm{k}', \\
    0, & \bm{k}\neq\bm{k}'.
    \end{array}
    \right.
\end{equation*}
Then, for a periodic function ${U}\in\mathcal{L}^2(\mathbb{T}^n)$, its Fourier series is defined by
\begin{equation*}\label{eq:parent Fourier series}
    {U}(\bm{y})=\sum_{\bm{k}\in\mathbb{Z}^n}\hat{{U}}_{\bm{k}}\varphi_{\bm{k}}(\bm{y}),\qquad
    \hat{{U}}_{\bm{k}}:=\left<{U},\varphi_{\bm{k}}\right>.
\end{equation*}

For a quasiperiodic function $u(\bm{x})\in\mathcal{Q}(\mathbb{R}^d)$, its mean value $M(u)$ is defined as 
$$M(u) := \lim_{T\rightarrow+\infty}\frac{1}{(2T)^d}\int_{\bm{s}+[-T,T]^d}u(\bm{x})\mathrm{d}\bm{x},\qquad \forall\bm{s}\in\mathbb{R}^d.$$
Correspondingly, the inner product and norm of $u,~v\in \mathcal{Q}(\mathbb{R}^d)$ can be defined as
\begin{equation*}\label{eq:inner norm}
    \begin{aligned}
    \left<u,v\right>:=M(u\bar{v}),\qquad \|u\|:=\left(M(|u|^2)\right)^{1/2}.\\
    \end{aligned}
\end{equation*}
We say $u\in\mathcal{L}^2_{\mathcal{Q}}$ if $\|u\|<+\infty$. 
Note that the definition of Fourier basis functions can be easily extended to more general cases as
$\varphi_{\bm{q}}(\bm{x}):=\mathrm{e}^{i\bm{q\cdot x}}$, for any $\bm{q},~\bm{x}\in\mathbb{R}^d$. 
Then, the Fourier-Bohr transform of $u$ is
$\hat{u}_{\bm{q}}:=M\big(u\varphi_{\bm{q}}\big),~\bm{q}\in\mathbb{R}^d$. 
\begin{lemma}[\cite{jiang2024numerical}, Theorem 4.1]\label{lem:general_Fourier}
    For a $d$-dimensional quasiperiodic function $u$ and its associated parent function $U$, it holds $\hat{u}_{\bm{q}}=\hat{{U}}_{\bm{k}}$ when $\bm{q}=\bm{Pk}$.
\end{lemma}
Then the generalized Fourier series of $u(\bm{x})\in \mathcal{L}^2_{\mathcal{Q}}$ is given by
\begin{equation*}\label{eq:QP Fourier series}
    u(\bm{x})=\sum_{\bm{k}\in\mathbb{Z}^n}\hat{{U}}_{\bm{k}}\varphi_{\bm{Pk}}(\bm{x}).
\end{equation*}

\begin{remark}\label{eq:fourier_decay}
    When the parent function satisfies certain regularity condition, all Fourier coefficients have the decay property.
    Specifically, if $U\in \mathcal{H}^{\alpha}(\mathbb{T}^n)$, there exists a positive constant $C$ such that $\hat{U}_{\bm{k}}\leq C|\bm{k}|^{-\alpha}|U|_{\mathcal{H}^{\alpha}}$~\cite{grafakos2008classical}. 
    Here, the definitions of seminorm $|\cdot|_{\mathcal{H}^{\alpha}}$ and the corresponding Sobolev space $\mathcal{H}^{\alpha}(\mathbb{T}^n)$ are given in \cref{eq:periodic norm}.
\end{remark}

Next, we briefly introduce the PM.
Unlike previous numerical methods, PM grasps the essential feature of a $d$-dimensional quasiperiodic function that can be embedded into its associated $n$-dimensional parent periodic function\,\cite{jiang2014numerical}. 
As a result, PM computes the $n$-dimensional parent periodic system in a pseudospectral way instead of directly addressing quasiperiodic system. Then, PM projects these results onto $d$-dimensional space by the projection matrix $\bm{P}$ to obtain quasiperiodic system. 
Concretely, given a positive integer $N$, we define the finite index set

$$
    \mathcal{K}_N:=\left\{\bm{k}\in\mathbb{Z}^n:\bm{k}\in[-N,N)^n\right\}.
$$
Then, the dual grid of $\mathcal{K}_N$ is given by
\begin{equation*}
    \mathcal{G}_N:=\left\{\bm{y_{\ell}}= \pi\bm{\ell}/N \in[0,2\pi)^n:\bm{\ell}\in\mathbb{Z}^n\cap[0,2N)^n\right\}.
\end{equation*}
For periodic functions $U$ and $V$, the compound trapezoidal formula of inner product is 
\begin{equation*}
    \left< U,V\right>_N:=\frac{1}{(2N)^n}\sum_{\bm{y_\ell}\in\mathcal{G}_N}U(\bm{y_\ell})\overline{V(\bm{y_\ell})}.
\end{equation*}
Limiting the space $\mathcal{L}^2(\mathbb{T}^n)$ to a finite dimensional subspace spanned by the $\{\varphi_{\bm{k}}:\bm{k}\in\mathcal{K}_N\}$, we obtain the discrete Fourier-Bohr series of $u\in \mathcal{L}^2_{\mathcal{Q}}$~\cite{jiang2024numerical}
\begin{equation*}
    u(\bm{x})= \sum_{\bm{k}\in\mathcal{K}_{N}}\bar{{U}}_{\bm{k}}\varphi_{\bm{Pk}}(\bm{x}),\qquad \bar{{U}}_{\bm{k}}:=\left<{U},\varphi_{\bm{k}}\right>_{N}.
\end{equation*}

The corresponding error analysis of PM can refer to \cite{jiang2024numerical}. Moreover, since discrete Fourier coefficients originate from the periodic parent function, PM can use the $n$-dimensional FFT to improve the computational efficiency.

\section{Irrational-window-filter projection method (IWFPM)}\label{sec:IWFPM}

When using PM to address some quasiperiodic systems, such as QSEs, an interesting phenomenon has been observed that the Fourier coefficients are concentrated in a narrow elongated area (see \cite{wang2022convergence,gao2023reduced}). 
In this section, based on this phenomenon, we improve the index set of spectral points, thereby reducing the DOF of PM. 
Then we further overcome the challenge of performing FFT on irregular index sets by using an index transform. 
Finally, we provide implementation details of IWFPM and establish an error analysis of this method at the functional level.

\subsection{Irrational window}\label{subsec:irr}

Based on the decay property of Fourier coefficients mentioned in \cref{eq:fourier_decay}, we discover that Fourier coefficients are concentrated in a hyperparallelogram area along the $\bm{Pk}=\bm{0}$ direction. 
Further, we divide $\bm{P}=(\bm{P}_{\mathrm{I}},\bm{P}_{\mathrm{II}})$, where $\bm{P}_{\mathrm{I}}\in\mathbb{R}^{d\times d}$ and $\bm{P}_{\mathrm{II}}\in\mathbb{R}^{d\times (n-d)}$.
According to the definition of projection matrix, we can always make the $d$-order matrix $\bm{P}_{\mathrm{I}}$ invertible.
Hence, there exists an elementary row transform $\bm{P}_{\mathrm{I}}^{-1}$ such that
$\bm{P}_{\mathrm{I}}^{-1}\bm{P}=(\mathbf{I}_d,\bm{Q})$, where $\mathbf{I}_d$ is the $d$-order identity matrix and $\bm{Q}:=\bm{P}_{\mathrm{I}}^{-1}\bm{P}_{\mathrm{II}}\in\mathbb{R}^{d\times(n-d)}$. 
Through this transform, we can concentrate all irrational numbers in $\bm{P}$ into $\bm{Q}$. 
Correspondingly, we partition the index $\bm{k}\in\mathbb{R}^n$ into two parts: $\bm{k}=(\bm{k}_{\mathrm{I}}^T,\bm{k}_{\mathrm{II}}^T)^T,~\bm{k}_{\mathrm{I}}\in\mathbb{R}^d,~\bm{k}_{\mathrm{II}}\in\mathbb{R}^{n-d}$. 
Then, the hyperparallelogram tilt along the $\bm{P}_{\mathrm{I}}^{-1}\bm{Pk}=\bm{k}_{\mathrm{I}}+\bm{Qk}_{\mathrm{II}}=\bm{0}$ direction. 

Based on this distribution feature of Fourier coefficients, we can define an irrational window for given two positive integers $K$ and $L$ as
\begin{equation*}
    \mathcal{W}_{K,L}:=\left\{\bm{k}=(\bm{k}_{\mathrm{I}}^T,\bm{k}_{\mathrm{II}}^T)^T\in\mathbb{R}^n:\bm{k}_{\mathrm{II}}\in[-L,L)^{n-d},~\bm{k}_{\mathrm{I}} +\bm{Q}\bm{k}_{\mathrm{II}} \in[-K,K)^d\right\}.
\end{equation*}
Obviously, irrational window $\mathcal{W}_{K,L}$ is determined by the irrational numbers in the projection matrix $\bm{P}$. 
Then, we can define a hyperparallelogram index set 
\begin{equation}\label{eq:non-tensor}
    \mathcal{K}_{K,L}:=\mathcal{W}_{K,L}\cap\mathbb{Z}^n.
\end{equation}

\begin{figure*}[!hbpt]
    \centering
    \includegraphics[width=0.45\textwidth]{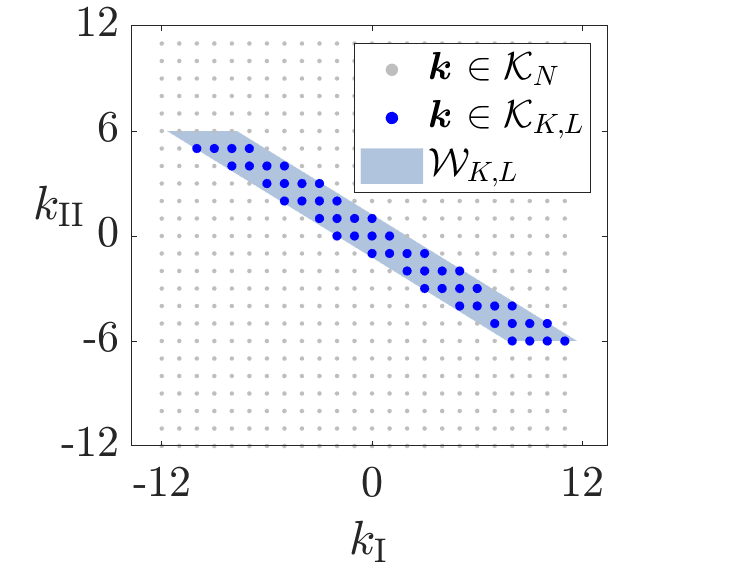}
    \caption{Rectangle index set $\mathcal{K}_N$, parallelogram index set  $\mathcal{K}_{K,L}$, and irrational window $\mathcal{W}_{K,L}$ when $d=1$, $n=2$, $N=12$, $K=2$, $L=6$, $\bm{P}=(1,(\sqrt{5}+1)/2)$.}
    \label{fig:window}
\end{figure*}

As an illustrative example, \Cref{fig:window} presents the rectangle index set $\mathcal{K}_N$, parallelogram index set $\mathcal{K}_{K,L}$, and irrational window $\mathcal{W}_{K,L}$ when $d=1$, $n=2$, $N=12$, $K=2$, $L=6$ and projection matrix $\bm{P}=(1,(\sqrt{5}+1)/2)$.
We can observe that the index set $\mathcal{K}_N$ has $2N\times2N=576$ points, while $\mathcal{K}_{K,L}$ significantly reduces the number of points to $2K\times2L=48$.

\begin{remark}
    Generally, the two positive integers $K$ and $L$ can be two vectors $\bm{K}=\{K_1,\cdots,K_d\}\in\mathbb{N}^d_+$ and $\bm{L}=\{L_1,\cdots,L_{n-d}\}\in\mathbb{N}^{n-d}_+$, respectively. 
\end{remark}

\subsection{Index-shift map}\label{subsec:index-filter}
For the Fourier coefficient index set $\mathcal{K}_{K,L}$, there seems a drawback in practical calculations that the irregular shape may make the FFT inapplicable. To address this issue, we introduce an index-shift map $\varrho$, defined by
\begin{equation*}\label{eq:index trans}
\varrho(\bm{k})=\bm{k}^*=(k^*_j)_{j=1}^n,
\end{equation*}
where
\begin{equation}\label{eq:index mod}
k^*_j=\left\{
\begin{array}{ll}
  k_j~{\rm mod}~2K, & \text{if}~ j=1,\cdots,d, \\[3mm]
  k_j~{\rm mod}~2L, & \text{if}~ j=d+1,\cdots,n.
\end{array}
\right.
\end{equation}
Here, ``${\rm mod}$" represents modulo operation. 

Applying $\varrho$ to all indicators of $\mathcal{K}_{K,L}$, we obtain the following hyperrectangle index set
\begin{equation*}\label{eq:hyperrectangle index set}
    \mathcal{K}^*_{K,L}:=\left\{\bm{k}=(\bm{k}_{\mathrm{I}}^T,\bm{k}_{\mathrm{II}}^T)^T\in\mathbb{Z}^n:\bm{k}_{\mathrm{I}}\in[0,2K)^d,~\bm{k}_{\mathrm{II}}\in[0,2L)^{n-d}\right\}.
\end{equation*}
Correspondingly, the set of dual grid points can be defined as
\begin{equation}\label{eq:G_KL}
\mathcal{G}_{K,L}:=\{\bm{y_{\ell}}=\left( \pi \bm{\ell}^T_1/K, \pi \bm{\ell}^T_2/L \right)^T\in[0,2\pi)^n:\bm{\ell}=(\bm{\ell}^T_1,\bm{\ell}^T_2)^T\in \mathcal{K}^*_{K,L}\}.
\end{equation}
Obviously, 
\begin{equation}\label{eq:basis_equivalent}
    \varphi_{\bm{k}}(\bm{y_l})=\varphi_{\bm{k}^*}(\bm{y_l}),\qquad\forall\bm{y_{\ell}}\in\mathcal{G}_{K,L},~\bm{k}\in\mathcal{K}_{K,L}.
\end{equation}

\begin{remark}
    Based on the equivalence relationship \cref{eq:basis_equivalent}, we can establish the connection between the discrete Fourier transforms on indicator set $\mathcal{K}_{K,L}$ and indicator set $\mathcal{K}^*_{K,L}$ in \Cref{subsec:details IWFPM}.
\end{remark}

\begin{figure*}[!hbpt]
    \centering
    \subfigure[$\mathcal{K}_{K,L}$]{\label{fig:non-tensor}\includegraphics[width=0.3\textwidth]{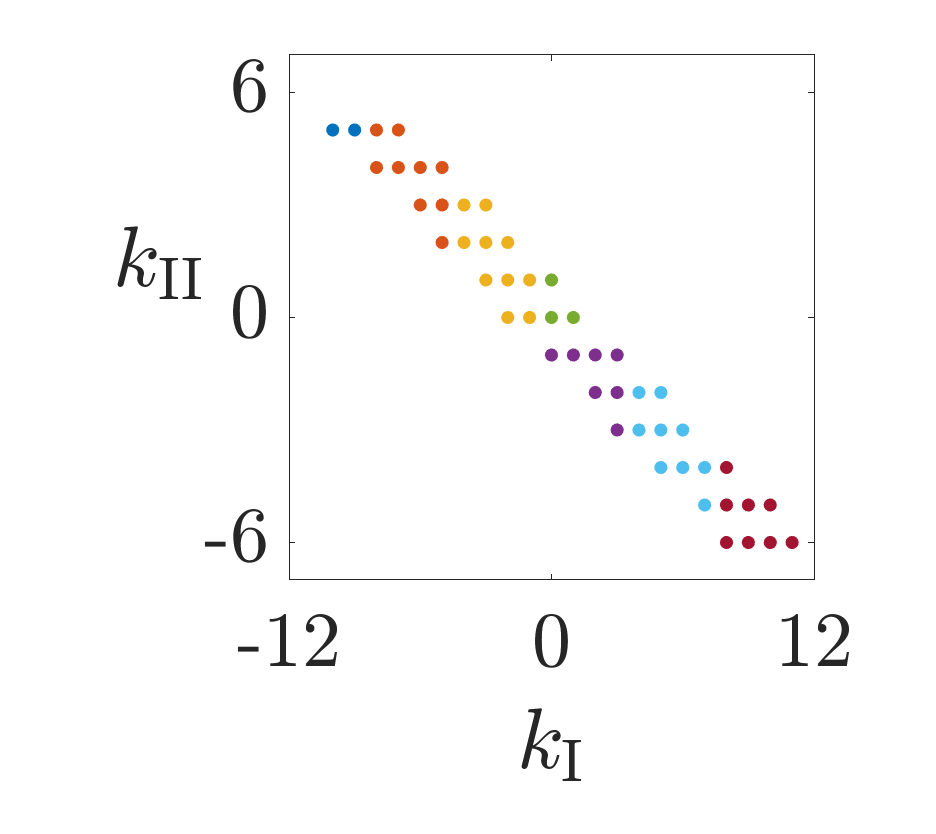}}
    \subfigure[$\mathcal{K}^*_{K,L}$]{\label{fig:tensor}\includegraphics[width=0.3\textwidth]{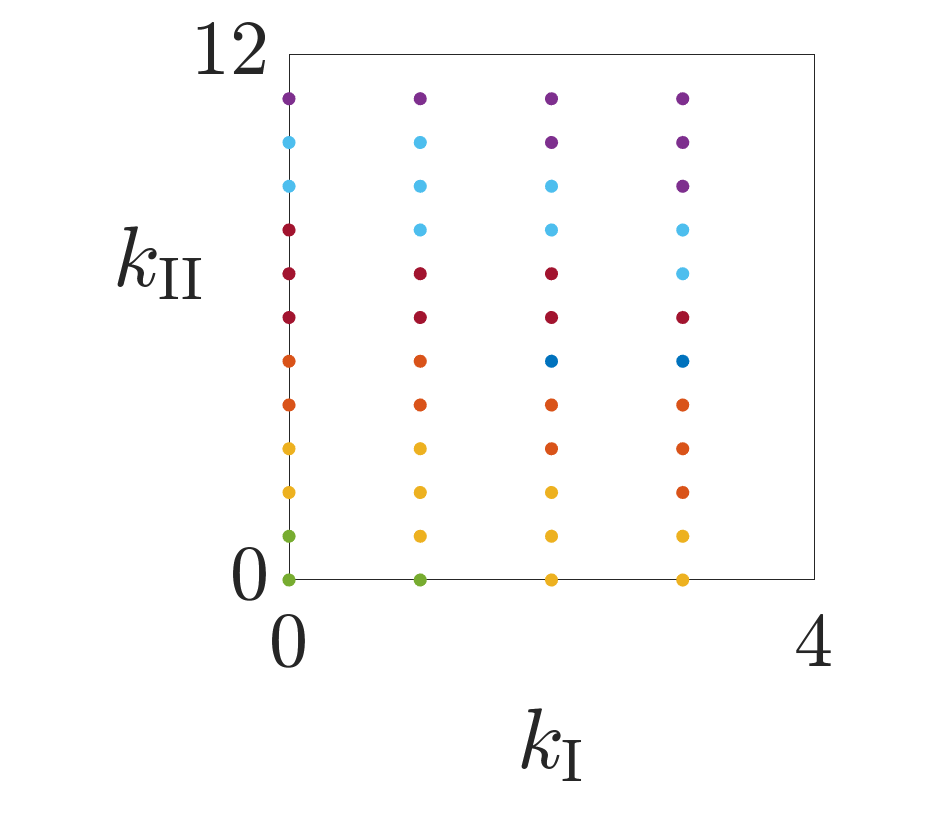}}
    \subfigure[$\mathcal{G}_{K,L}$]{\label{fig:grid}\includegraphics[width=0.3\textwidth]{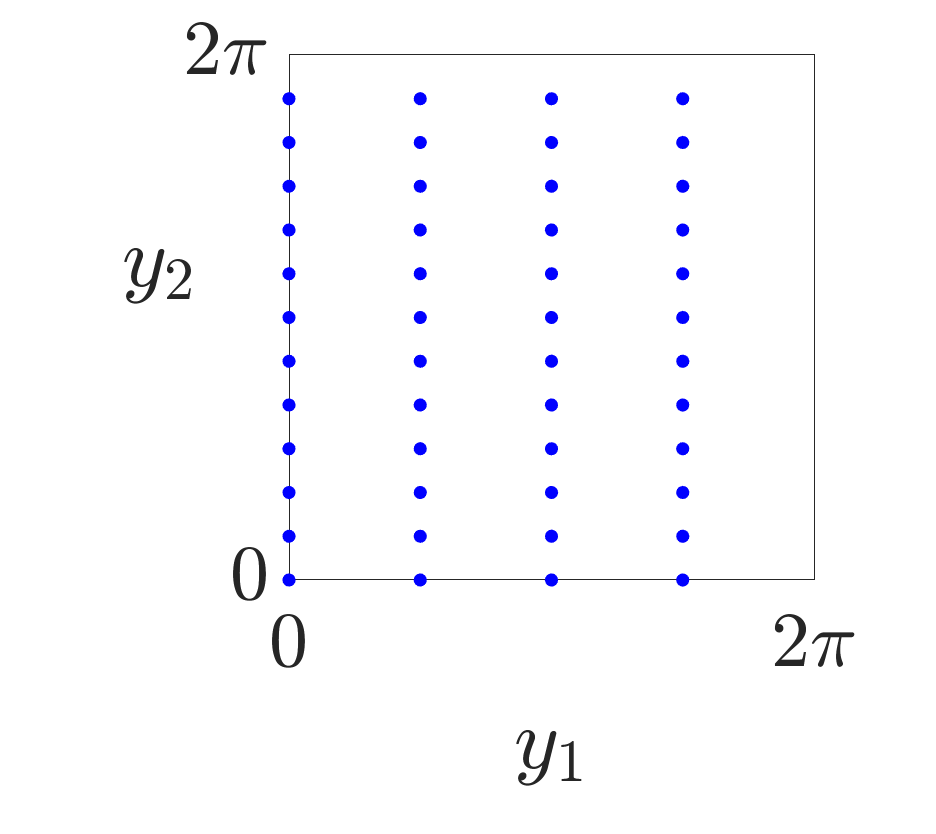}}
        \caption{The parallelogam index set $\mathcal{K}_{K,L}$ (left), the rectangle index set $\mathcal{K}^*_{K,L}$ (middle), and the grid points $\mathcal{G}_{K,L}$ (right) when $d=1$, $n=2$, $K=2$, $L=6$, $\bm{P}=(1,(\sqrt{5}+1)/2)$.}
    \label{fig:points}
\end{figure*}

As an example, \Cref{fig:points} illustrates the index-shift map on the index set $\mathcal{K}_{K,L}$ when $d=1$, $n=2$, $K=2$, $L=6$, and $\bm{P}=(1,(\sqrt{5}+1)/2)$. It shows the parallelogam index set $\mathcal{K}_{K,L}$, the rectangle index set $\mathcal{K}^*_{K,L}$, and the grid points $\mathcal{G}_{K,L}$. 

In what follows, we give an explicit expression of $\varrho^{-1}$ in another way. 
Although $\varrho$ is a bijection, the calculation of the inverse of $\varrho$ cannot be directly obtained by \cref{eq:index mod}. 
Hence, we attempt to give the mathematical expression of $\varrho^{-1}$ in the following.  
Let $(\bm{Qk}_{\mathrm{II}})_j$ denote the $j$th component of the vector $\bm{Qk}_{\mathrm{II}}$. 
According to the definition of $\bm{k}_{\mathrm{I}}=(k_j)_{j=1}^d$ in index $\bm{k}\in\mathcal{K}_{K,L}$, it is obvious that the value range of $k_j$ is the $2K$ integers within the interval $[-K-(\bm{Qk}_{\mathrm{II}})_j,K-(\bm{Qk}_{\mathrm{II}})_j)$. Among these integers, only one can be divisible by $2K$, and that is
$R_j=\lceil(-K-(\bm{Qk}_{\mathrm{II}})_j)/2K\rceil\cdot 2K$,
where $\lceil\cdot\rceil$ represents rounding up. 
Then, we can obtain the inverse map
\begin{equation}\label{eq:inverse1}
    \varrho^{-1}(\bm{k}^*)=\bm{k}=(\bm{k}_{\mathrm{I}}^T,\bm{k}_{\mathrm{II}}^T)^T,\qquad \bm{k}^*\in\mathcal{K}^*_{K,L}
\end{equation}
by two steps. The first step is to compute $\bm{k}_{\mathrm{II}}$ as following
\begin{equation}\label{eq:inverse2}
    k_j=\left\{
\begin{array}{ll}
  k_j^*, & \text{if}~k_j^*<L, \\
  k_j^*-2L, & \text{if}~L\leq k_j^*<2L,
    \end{array}
    \right.~j=d+1,\cdots,n,
\end{equation}
and the second step is to compute $\bm{k}_{\mathrm{I}}$ as following
\begin{equation}\label{eq:inverse3}
    k_j=\left\{
    \begin{array}{ll}
    k_j^*+R_j, & \text{if}~k_j^*+R_j+(\bm{Qk}_{\mathrm{II}})_j<K, \\
    k_j^*+R_j-2K, & \text{otherwise},
    \end{array}
\right.~ j=1,\cdots,d.
\end{equation}

\subsection{Implementation}\label{subsec:details IWFPM}
Based on the new index set $\mathcal{K}_{K,L}$ and the index-shift map $\varrho$ introduced in the previous two subsections, we can now present the implementation of IWFPM.

For $n$-dimensional periodic functions $U$ and $V$, the compound trapezoidal formula of inner product over $\mathcal{G}_{K,L}$ is
\begin{equation}\label{eq:new_inner_product}
    \begin{aligned}
    \left<{U},{V}\right>_{K,L}:=\frac{1}{2^nK^dL^{n-d}}\sum_{\bm{y_\ell}\in\mathcal{G}_{K,L}}{U}(\bm{y_\ell})\overline{{V}(\bm{y_\ell})}.
    \end{aligned}
\end{equation}
The finite dimensional subspace of $\mathcal{L}^2(\mathbb{T}^n)$ is
$
    \mathcal{S}_{K,L}:={\rm span}\{\varphi_{\bm{k}}:\bm{k}\in\mathcal{K}_{K,L}\}.
$
Then, the Fourier interpolation operator of IWFPM is defined by
\begin{equation}\label{eq:int}
\begin{aligned}
\mathcal{I}_{K,L}:~~\mathcal{L}^2_{\mathcal{Q}}&\rightarrow \mathcal{S}_{K,L},\\
u(\bm{x})&\mapsto \sum_{\bm{k}\in\mathcal{K}_{K,L}}\tilde{{U}}_{\bm{k}}\varphi_{\bm{Pk}}(\bm{x}),
\end{aligned}
\end{equation}
where $\tilde{{U}}_{\bm{k}}:=\left<{U},\varphi_{\bm{k}}\right>_{K,L}
=\left<{U},\varphi_{\bm{k}^*}\right>_{K,L}
$. 

Denoting $\mathbf{U}:=\left({U}(\bm{y_\ell})\right)_{\bm{y_\ell}\in\mathcal{G}_{K,L}}$,
there is a discrete Fourier transform $\mathbf{F}$ equivalent to $\mathcal{I}_{K,L}$ such that
\begin{equation}\label{eq:F}
    \mathbf{\tilde{U}}=\mathbf{F}\mathbf{U},\qquad \mathbf{\tilde{U}}:=\big(\tilde{{U}}_{\bm{k}}\big)_{\bm{k}\in\mathcal{K}_{K,L}}.
\end{equation}
Correspondingly, let $\mathbf{F}^*$ be the standard discrete Fourier transform, i.e. 
\begin{equation}\label{eq:standard discrete Fourier}
    \mathbf{\tilde{U}^*}=\mathbf{F}^*\mathbf{U},\qquad\mathbf{\tilde{U}^*}:=\big(\tilde{{U}}^*_{\bm{k}^*}\big)_{\bm{k}^*\in\mathcal{K}^*_{K,L}},
\end{equation}
where $\tilde{{U}}^*_{\bm{k}^*}:=\left<{U},\varphi_{\bm{k}^*}\right>_{K,L}$. 
Based on \Cref{subsec:index-filter}, 
we can define the index-shift transform $\mathbf{T}$ with respect to the inverse map $\varrho^{-1}$ such that
$
    \mathbf{\tilde{U}}=\mathbf{T}\mathbf{\tilde{U}^*}.
$ 
Therefore, to implement the discrete Fourier transform $\mathbf{F}$, we can equivalently apply the transform $\mathbf{F}^*$ and $\mathbf{T}$ successively. 
In other words, it means that
$
\mathbf{F}=\mathbf{TF}^*
$. 
Through the index-shift transform $\mathbf{T}$, 
FFT is available to be performed on index set $\mathcal{K}^*_{K,L}$ with DOF $(2K)^d(2L)^{n-d}$.

\Cref{alg:IWFPM} summarizes the implementation process of IWFPM.
\begin{algorithm}[!pbht]
	\caption{Irrational-window-filter projection method (IWFPM)}
	\label{alg:IWFPM}
	\begin{algorithmic}[1]
		\REQUIRE projection matrix $\bm{P}$, size of index set: $K$ and $L$
		\STATE Generate index sets $\mathcal{K}_{K,L}$ and $\mathcal{K}^*_{K,L}$
		\STATE Obtain  $\mathbf{\tilde{U}^*}=\big(\tilde{{U}}^*_{\bm{k}^*}\big)_{\bm{k}^*\in\mathcal{K}^*_{K,L}}$ by solving \cref{eq:standard discrete Fourier} using FFT
		\FOR{$\bm{k}^*\in\mathcal{K}^*_{K,L}$}
		\STATE Solve $\bm{k}=\varrho^{-1}(\bm{k}^*)$ according to \cref{eq:inverse1}-\cref{eq:inverse3}
		\STATE Store $\tilde{{U}}_{\bm{k}}=\tilde{{U}}^*_{\bm{k}^*}$
		\ENDFOR
		\STATE Calculate the Fourier interpolation $u=\sum_{\bm{k}\in\mathcal{K}_{K,L}}\tilde{{U}}_{\bm{k}}\varphi_{\bm{Pk}}$
	\end{algorithmic}
\end{algorithm}

\begin{remark}
    Note that the index-shift operator $\mathbf{T}$ only modifies the indicators of Fourier coefficients without changing the value. 
    Thus, the computational cost of the discrete Fourier transform $\mathbf{F}$  corresponding to IWFPM is entirely equivalent to that of the standard discrete Fourier transform $\mathbf{F}^*$. 
    Further, by utilizing FFT, the computational cost of $\mathbf{F}$ is of order $O\big(K^dL^{n-d}(\log K+\log L)\big).$ 
\end{remark}

\subsection{Error analysis}\label{subsec:error}

In this subsection, we give an error analysis of IWFPM. 
For simplicity of analysis, we use the notation $A[u]\lesssim B[u]$,
which means that there exists a positive constant satisfying $A[u]\leq CB[u]$, where $A[u]$ and $B[u]$ are functional with respect to $u(\bm{x})$, and the positive constant $C$ is independent of $K$, $L$, and any norm of $u$. 
Moreover, we abbreviate $A[u]\lesssim B[u]$ and $B[u]\lesssim A[u]$ as $A[u]\simeq B[u]$. 
The Sobolev seminorm and norm of quasiperiodic function $u\in\mathcal{L}^2_\mathcal{Q}$ for any $\alpha\geq0$ are defined as
\begin{equation}\label{eq:QP Sobolev norm}
\begin{aligned}
|u|_{\alpha}&:=\bigg(\sum_{\bm{k}\in\mathbb{Z}^n}\|\bm{Pk}\|^{2\alpha}|\hat{{U}}_{\bm{k}}|^2\bigg)^{1/2},\\
\|u\|_{\alpha}&:=\bigg(\sum_{\bm{k}\in\mathbb{Z}^n}(1+\|\bm{Pk}\|^{2\alpha})|\hat{{U}}_{\bm{k}}|^2\bigg)^{1/2}.
\end{aligned}
\end{equation}
Here we set $0^0=1$. And we say $u\in\mathcal{H}^\alpha_{\mathcal{Q}}$ if $\|u\|_{\alpha}<+\infty$. 
The Sobolev seminorm and norm of periodic function ${U}\in\mathcal{L}^2(\mathbb{T}^n)$ for any $\alpha\geq0$ are defined as
\begin{equation}\label{eq:periodic norm}
\begin{aligned}
|{U}|_{\mathcal{H}^\alpha}&:=\bigg(\sum_{\bm{k}\in\mathbb{Z}^n}\|\bm{k}\|^{2\alpha}|\hat{{U}}_{\bm{k}}|^2 \bigg)^{1/2},\\
\|{U}\|_{\mathcal{H}^\alpha}&:=\bigg(\sum_{\bm{k}\in\mathbb{Z}^n}\big(1+\|\bm{k}\|^{2\alpha}\big)|\hat{{U}}_{\bm{k}}|^2\bigg)^{1/2}. 
\end{aligned}
\end{equation}
And we say ${U}\in\mathcal{H}^\alpha(\mathbb{T}^n)$ if $\|{U}\|_{\mathcal{H}^\alpha}<+\infty$. 
If the norm of the projection matrix $\bm{P}$ is not very large, the quasiperiodic norm $\|u\|_{\alpha}$ can be effectively controlled by the periodic norm $\|{U}\|_{\mathcal{H}^\alpha}$, while the opposite is not true. $u\in\mathcal{H}^\alpha_{\mathcal{Q}}$ may not necessarily lead to ${U}\in\mathcal{H}^\alpha(\mathbb{T}^n)$, or the norm $\|{U}\|_{\mathcal{H}^\alpha}$ may be much larger than the norm $\|u\|_{\alpha}$. Considering the definition of index set $\mathcal{K}_{K,L}$ \cref{eq:non-tensor}, we adopt a new norm definition
\begin{equation}\label{eq:new norm}
\begin{aligned}
|u|_{\alpha,\beta}&:=\bigg(\sum_{\bm{k}\in\mathbb{Z}^n}\big(\|\bm{k}_{\mathrm{I}} +\bm{Qk}_{\mathrm{II}}\|^{2\alpha}+ \|\bm{k}_{\mathrm{II}}\|^{2\beta}\big)|\hat{{U}}_{\bm{k}}|^2\bigg)^{1/2},\\
\|u\|_{\alpha,\beta}&:=\bigg(\sum_{\bm{k}\in\mathbb{Z}^n}\big(1+\|\bm{k}_{\mathrm{I}}+\bm{Qk}_{\mathrm{II}}\|^{2\alpha}+ \|\bm{k}_{\mathrm{II}}\|^{2\beta}\big)|\hat{{U}}_{\bm{k}}|^2\bigg)^{1/2},
\end{aligned}
\end{equation}
for any $\alpha,\beta\geq0$. Note that the Cauchy-Schwarz inequality can lead to
\begin{footnotesize}
$$\begin{aligned}
\bigg(\sum_{\bm{k}\in\mathbb{Z}^n}\lambda_{\bm{k}}^2|\hat{U}_{\bm{k}}+\hat{V}_{\bm{k}}|^2\bigg)^{1/2}
\leq&\bigg(\sum_{\bm{k}\in\mathbb{Z}^n}\lambda_{\bm{k}}^2\big(|\hat{U}_{\bm{k}}|^2+ |\hat{V}_{\bm{k}}|^2\big)+2\sum_{\bm{k}\in\mathbb{Z}^n}\lambda_{\bm{k}}^2 |\hat{U}_{\bm{k}}\hat{V}_{\bm{k}}|\bigg)^{1/2}\\
\leq&\bigg(\sum_{\bm{k}\in\mathbb{Z}^n}\lambda_{\bm{k}}^2\big(|\hat{U}_{\bm{k}}|^2+ |\hat{V}_{\bm{k}}|^2\big)+2\bigg(\sum_{\bm{k}\in\mathbb{Z}^n}\lambda_{\bm{k}}^2 |\hat{U}_{\bm{k}}|^2\sum_{\bm{k}\in\mathbb{Z}^n}\lambda_{\bm{k}}^2|\hat{V}_{\bm{k}}|^2\bigg)^{1/2}\bigg)^{1/2}\\
=&\bigg(\sum_{\bm{k}\in\mathbb{Z}^n}\lambda_{\bm{k}}^2 |\hat{U}_{\bm{k}}|^2\bigg)^{1/2}+ \bigg(\sum_{\bm{k}\in\mathbb{Z}^n}\lambda_{\bm{k}}^2 |\hat{V}_{\bm{k}}|^2\bigg)^{1/2}
\end{aligned}$$
\end{footnotesize}
for any real sequence $\{\lambda_{\bm{k}}\}_{\bm{k}\in\mathbb{Z}^n}$. Therefore, it is easy to prove that $|\cdot|_{\alpha,\beta}$ and $\|u\|_{\alpha,\beta}$ satisfy the conditions for defining the seminorms and the norm, respectively. Here we say $u\in\mathcal{H}^{\alpha,\beta}_{\mathcal{Q}}$ if $\|u\|_{\alpha,\beta}<+\infty$. 

Then we give the following lemma shows the relation among the above three spaces.

\begin{lemma}\label{lm:norm relation}
Suppose that $\alpha\geq\beta\geq0$, then $u\in\mathcal{H}^{\alpha,\beta}_{\mathcal{Q}}$ if and only if $u\in\mathcal{H}^\alpha_{\mathcal{Q}}$ and ${U}\in\mathcal{H}^\beta(\mathbb{T}^n)$, i.e. the seminorms and norms defined in \cref{eq:QP Sobolev norm}, \cref{eq:periodic norm} and \cref{eq:new norm} satisfy
$$\begin{aligned}
|u|_{\alpha,\beta}&\simeq |u|_{\alpha}+|{U}|_{\mathcal{H}^\beta},\\
\|u\|_{\alpha,\beta}&\simeq \|u\|_{\alpha}+\|{U}\|_{\mathcal{H}^\beta}.
\end{aligned}$$
\end{lemma}
\begin{proof}
The proof is in \Cref{sec:proof_norm_relation}.
\end{proof}

The truncation approximation operator can be defined via the following projection operator
\begin{equation}\label{eq:proj}
\begin{aligned}
\mathcal{P}_{K,L}~:~\mathcal{L}^2_{\mathcal{Q}}&\rightarrow \mathcal{S}_{K,L}\\
u&\mapsto \sum_{\bm{k}\in\mathcal{K}_{K,L}}\hat{{U}}_{\bm{k}}\varphi_{\bm{Pk}},
\end{aligned}
\end{equation}
where $\hat{{U}}_{\bm{k}}$ is the Fourier coefficient for $\bm{k}\in\mathcal{K}_{K,L}$.

\begin{lemma}\label{th:proj}
Suppose that $u\in\mathcal{H}^{\alpha,\beta}_{\mathcal{Q}}$ with $\alpha\geq\beta\geq0$, then the error of truncation approximation $\mathcal{P}_{K,L}$ \cref{eq:proj} satisfies
$$\begin{aligned}
\left|u-\mathcal{P}_{K,L}u\right|_{\mu,\nu}&\lesssim K^{-\alpha}\big(K^{\mu}+L^{\nu}\big)|u|_{\alpha}+L^{-\beta}\big(K^{\mu}+L^{\nu}\big)|{U}|_{\mathcal{H}^\beta},\\
\left\|u-\mathcal{P}_{K,L}u\right\|_{\mu,\nu}&\lesssim K^{-\alpha}\big(K^{\mu}+L^{\nu}\big)\|u\|_{\alpha}+L^{-\beta}\big(K^{\mu}+L^{\nu}\big)\|{U}\|_{\mathcal{H}^\beta},
\end{aligned}$$
for $\mu\in[0,\alpha]$, $\nu\in[0,\beta]$.
\end{lemma}
\begin{proof}
The proof is in \Cref{sec:proof_proj}.
\end{proof}

Combined \Cref{th:proj} with \Cref{lm:norm relation}, the following corollary can be easily obtained.
\begin{corollary}\label{cr:proj}
Suppose that $u\in\mathcal{H}^{\alpha,\beta}_{\mathcal{Q}}$ with $\alpha\geq\beta\geq0$, then the error of truncation approximation $\mathcal{P}_{K,L}$ \cref{eq:proj} satisfies
$$\begin{aligned}
\left|u-\mathcal{P}_{K,L}u\right|_{\mu,\nu}&\lesssim \big(K^{-\alpha}+L^{-\beta}\big)\big(K^{\mu}+L^{\nu}\big)|u|_{\alpha,\beta},\\
\left\|u-\mathcal{P}_{K,L}u\right\|_{\mu,\nu}&\lesssim \big(K^{-\alpha}+L^{-\beta}\big)\big(K^{\mu}+L^{\nu}\big)\|u\|_{\alpha,\beta},
\end{aligned}$$
for $\mu\in[0,\alpha]$, $\nu\in[0,\beta]$.
\end{corollary}

\begin{theorem}\label{th:int}
Suppose that $u\in\mathcal{H}^{\alpha,\beta}_{\mathcal{Q}}$ with $\alpha\geq\beta> \dfrac{n-d}{2}$ and $\dfrac{d}{2\alpha}+\dfrac{n-d}{2\beta}<1$,
then the error of interpolation approximation $\mathcal{I}_{K,L}$ \cref{eq:int} is 
$$\begin{aligned}
|u-\mathcal{I}_{K,L}u|_{\mu,\nu}&\lesssim \big(K^{-\alpha}+L^{-\beta}\big)\big(K^{\mu}+L^{\nu}\big)|u|_{\alpha,\beta},\\
\|u-\mathcal{I}_{K,L}u\|_{\mu,\nu}&\lesssim \big(K^{-\alpha}+L^{-\beta}\big)\big(K^{\mu}+L^{\nu}\big)\|u\|_{\alpha,\beta},
\end{aligned}$$
for $\mu\in[0,\alpha]$, $\nu\in[0,\beta]$.
\end{theorem}
\begin{proof}
The proof is in \Cref{sec:proof_int}.
\end{proof}

\begin{remark}
According to the above convergence result of IWFPM interpolation, it is evident that for some special quasiperiodic functions with distinct regularities along different directions (means that the gap between $\alpha$ and $\beta$ is huge), our proposed hyperparallelogam index set $\mathcal{K}_{K,L}$ can achieve the consistent convergence effect by adjusting $K$ and $L$.
\end{remark}

\section{Application to quasiperiodic Schr\"{o}dinger eigenproblems (QSEs)}\label{sec:num}

In this section, we apply IWFPM to solve 1D, 2D, and 3D QSEs. Considering the eigenproblems with quasiperiodic Schr\"{o}dinger operator  $H:\mathcal{C}^2(\mathbb{R}^d)\rightarrow\mathcal{C}(\mathbb{R}^d)$ as
\begin{equation}\label{eq:QSE}
    Hu(\bm{x}):=-\frac{1}{2}\Delta u(\bm{x})+v(\bm{x})u(\bm{x})=Eu(\bm{x}),
\end{equation}
where $v(\bm{x})$ is a quasiperiodic potential, the eigenfunction $u(\bm{x})$ is the normalized wavefunction, and the eigenvalue $E$ represents the corresponding energy. 

In terms of the existence of QSE solutions, considerable progress has been made for 1D operators on both $\mathbb{Z}$ and $\mathbb{R}$~\cite{simon2000schrodinger,avila2009ten,avila2015global,avila2017sharp,ge2023multiplicative}. 
While it becomes significantly difficult when dealing with multi-dimensional QSEs and only a few papers exist~\cite{bourgain2007anderson,shi2021absence}. 
Moreover, the regularity analysis of QSE solutions remains an open problem.
Due to the challenges in developing theoretical research on QSEs, the IWFPM method holds great significance, as it offers a way to predict the shape of the solution from a numerical perspective.

\subsection{IWFPM discretization}\label{subsec:num_imp}

Suppose that $U(\bm{y})$ and $V(\bm{y})$ are the parent functions of $u(\bm{x})$ and $v(\bm{x})$, respectively.
Let $\mathbf{\tilde{U}}$ be the Fourier coefficients vector of $U(\bm{y})$ on the index set $\mathcal{K}_{K,L}$. 
Then, by the discretization of IWFPM, solving QSE \cref{eq:QSE} can be expressed as finding an eigenpair $(E,\mathbf{\tilde{U}})$ such that
\begin{equation*}\label{eq:pseudospectral vector}
\mathbf{\tilde{H}}\mathbf{\tilde{U}}:=\bm{\Lambda}\mathbf{\tilde{U}}+\mathbf{F}\mathbf{V}\mathbf{F}^{-1}\mathbf{\tilde{U}}=E\mathbf{\tilde{U}}
\end{equation*}
where 
\begin{equation*}\label{eq:discrete_Lambda}
    \bm{\Lambda}=\frac12\diag\left(\|\bm{Pk}\|^2\right)_{\bm{k}\in\mathcal{K}_{K,L}},\quad
    \mathbf{V}=(V(\bm{y}))_{\bm{y}\in\mathcal{G}_{K,L}},
\end{equation*}
$\mathbf{F}$ is the discrete Fourier transform corresponding to IWFPM defined by \cref{eq:F}. 
To solve this eigenvalue problem in matrix form, we employ the locally optimal block preconditioned conjugate gradient (LOBPCG) method\,\cite{knyazev2001toward}, with  convergence error 1.0e-10 and initial vector $\bm{e}_1=(1,0,\cdots,0)^T$. 
The preconditioner selected in LOBPCG method is 
\begin{equation}\label{eq:precond}
    \bm{M}=\argmin_{\bm{D}\in\mathcal{D}}\|\mathbf{\tilde{H}}\bm{D}-\bm{I}\|_F=
    \diag\big(   \tilde{h}_{11}/\|\mathbf{\tilde{H}}\bm{e}_1\|^2_2, 
    \cdots,
    \tilde{h}_{NN}/\|\mathbf{\tilde{H}}\bm{e}_N\|^2_2
   \big),
\end{equation}
where $N=(2K)^d(2L)^{n-d}$ is the size of matrix $\mathbf{\tilde{H}}$, and $\tilde{h}_{ii}$ is the $i$-th diagonal element of $\mathbf{\tilde{H}}$, $i=1,\cdots,N$. $\|\cdot\|_F$ means the Frobenius norm and $\mathcal{D}$ is the set of all diagonal matrices of order $N$. More details about this preconditioner can refer to \cite{jiang2024high}. 
\cref{alg:IWFPM for QSE} summarizes the detailed process of using LOBPCG eigensolver to solve the QSE \cref{eq:QSE}. 
Note that, the LOBPCG eigensolver can simultaneously compute multiple eigenpairs of arbitrary QSE.

\begin{algorithm}[!pbht]
	\caption{LOBPCG eigensolver for solving QSE \cref{eq:QSE}}
	\label{alg:IWFPM for QSE}
	\begin{algorithmic}[1]
		\REQUIRE $\bm{\Lambda}$, $\mathbf{V}$, $\bm{M}$, an initial vector $\mathbf{\tilde{U}}^{(0)}$, a conjugate vector $\mathbf{p}^{(0)}=0$
		\STATE Itetate: For $i=0,\cdots,$ until convergence:
		\STATE $\quad\mu^{(i)}:=\left< \mathbf{\tilde{U}}^{(i)},\mathbf{\tilde{U}}^{(i)}\right>/\left< \mathbf{\tilde{U}}^{(i)},\mathbf{\tilde{H}}\mathbf{\tilde{U}}^{(i)}\right>$, where $\mathbf{\tilde{H}}\mathbf{\tilde{U}}^{(i)}=\bm{\Lambda}\mathbf{\tilde{U}}^{(i)}+\mathbf{F}\mathbf{V}\mathbf{F}^{-1}\mathbf{\tilde{U}}^{(i)}$
		\STATE $\quad\mathbf{r}:=\mathbf{\tilde{U}}-\mu^{(i)}\mathbf{\tilde{H}}\mathbf{\tilde{U}}^{(i)}$
		\STATE $\quad\mathbf{w}=\bm{M}\mathbf{r}$
		\STATE $\quad$Use the Rayleigh-Ritz method for $\mathbf{I}-\mu^{(i)}\mathbf{\tilde{H}}$ on the trial subspace Span $\left\{\mathbf{w}^{(i)},\mathbf{\tilde{U}}^{(i)},\mathbf{p}^{(i)}\right\}$
		\STATE $\quad\mathbf{\tilde{U}}^{(i+1)}:=\mathbf{w}^{(i)}+\tau^{(i)}\mathbf{\tilde{U}}^{(i)}+\gamma^{(i)}\mathbf{p}^{(i)}$(the Ritz vector corresponding to the maximal Ritz value)
		\STATE $\quad\mathbf{p}^{(i+1)}:=\mathbf{w}^{(i)}+\gamma^{(i)}\mathbf{p}^{(i)}$
		\STATE End 
		\STATE Output the approximations $E=1/\mu^{(i)}$ and $\mathbf{\tilde{U}}=\mathbf{\tilde{U}}^{(i)}$ to the smallest eigenvalue and its corresponding eigenvector.
	\end{algorithmic}
\end{algorithm}

\subsection{Numerical experiments}\label{subsec:num_exp}
Now we present the numerical results obtained by IWFPM and demonstrate the performance by comparing with PM. 
All algorithms are coded by Matlab 2022b. The computations for 1D and 2D QSEs are carried out on a workstation with an Inter Core 2.10 GHz CPU and 16 GB RAM. IT and CPU represent the required iterations and the CPU time (in seconds), respectively. $\text{DOF}:=(2K)^d(2L)^{n-d}$ denotes the degrees of freedom. 
In this section, we consistently present the calculation results of the minimum eigenvalue $E_0$ and the corresponding eigenfunction $u_0(\bm{x})$. 
We observe the eigfunction in grid form as $\bm{u}_0=(u_0(\bm{\xi}))_{\bm{\xi}\in\mathcal{G}}$, where $\mathcal{G}$ is a uniform grid on the bounded region $\mathcal{G}=[-5a,5a]^d~(a=10^{3-d})$ with the step size $h=0.1$, 
and normalize it through dividing by the norm of maximal module $\|\bm{u}_0\|_{\infty}$.
The probability density of eigenfunction $\bm{u}_0$ is denoted as $\bm{\rho}:=|\bm{u}_0|^2$. 
We use relative errors of $E_0$ and $\bm{u}_0$
to measure the numerical accuracy
\begin{equation*}
    E_v=\bigg|\frac{E_0-E_0^*}{E_0^*}\bigg| 
    \quad \text{and}\quad 
    E_f=\|\bm{u}_0-\bm{u}_0^*\|_\infty,
\end{equation*}
where $E^*_0$ and $\bm{u}^*_0$ are corresponding numerical exact solutions of $E_0$ and $\bm{u}_0$, respectively. 

\begin{example}\label{exm:1}
    Consider 1D QSE \cref{eq:QSE} with potential
    \begin{equation}\label{eq:E1_equation}
        v(x)=v_0[2-\cos(2\pi x)-\cos(2\pi\alpha x)],
    \end{equation}
    where $v_0\in\mathbb{R},~\alpha=(\sqrt{5}-1)/2$. 
\end{example}
The projection matrix corresponding to \cref{eq:E1_equation} is $\bm{P}=2\pi(1,\alpha)$, then $v(x)$ can be embedded into the 2D parent function $V(\bm{y})=v_0(2-\cos{y_1}-\cos{y_2}),\,\bm{y}=(y_1,y_2)^T$. According to \cref{eq:non-tensor}, the hyperparallelogram index set is
\begin{equation*}
    \mathcal{K}_{K,L}=\left\{\bm{k}=(k_1,k_2)^T\in\mathbb{Z}^2:k_1+\alpha k_2\in[-K,K),~k_2\in[-L,L)\right\}.
\end{equation*}

\begin{figure*}[!hbpt]
    \centering
    \subfigure[$v_0=2.5$]{\label{fig:E1_wave1}\includegraphics[width=0.35\textwidth]{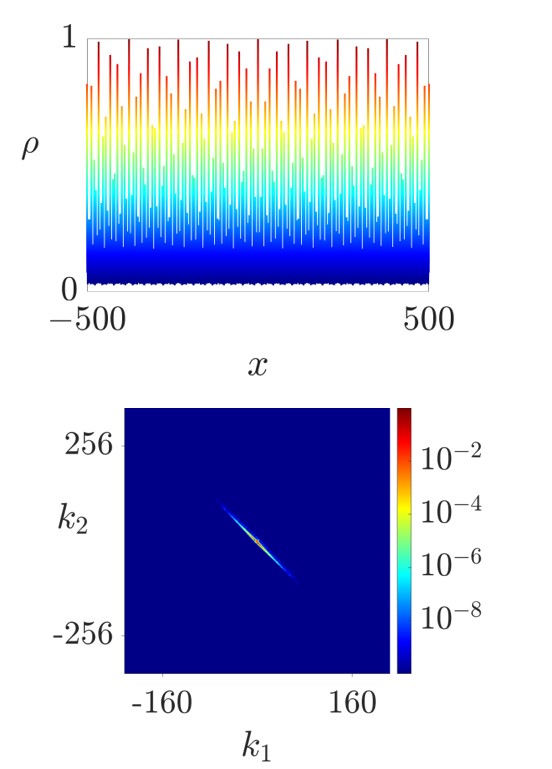}}
    \subfigure[$v_0=3$]{\label{fig:E1_wave2}\includegraphics[width=0.34\textwidth]{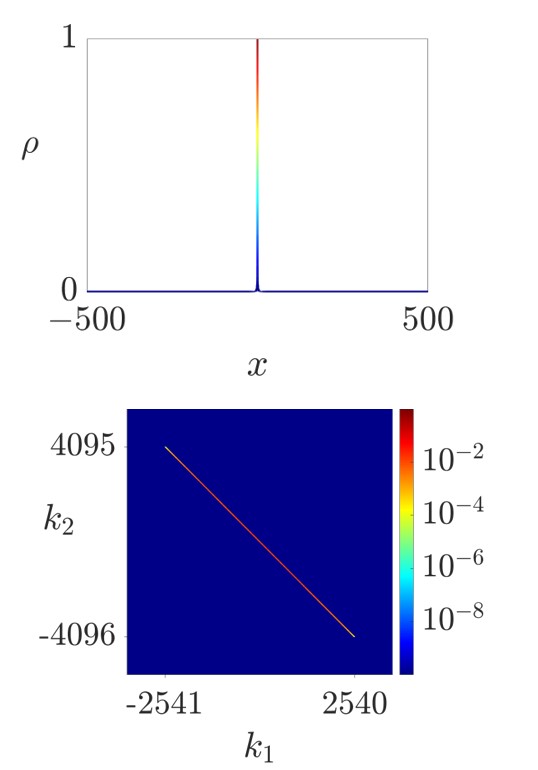}}
    \caption{Results of solving 1D QSE with potential \cref{eq:E1_equation} by IWFPM. The top row: probability density function $\bm{\rho}$; The bottom row: Fourier coefficients $\tilde{U}_{\bm{k}}$.}
    \label{fig:E1_wave}
\end{figure*}

This example is worth considering due to the observable phase transition from extended state to localized state as $v_0$ increases.  
To verify this, we present the probability density function $\bm{\rho}$ and the generalized Fourier coefficients $\tilde{U}_{\bm{k}}$ under the potential \cref{eq:E1_equation} with different $v_0$, as shown in \Cref{fig:E1_wave}. 
The wavefunction exhibits an extended state when $v_0=2.5$, and translates into a localized state when $v_0=3$. 
Moreover, it can be observed that the Fourier coefficients $\tilde{U}_{\bm{k}}$, whose intensities are larger than 1.0e-8, are mainly concentrated within a narrow parallelogram area. 
IWFPM method has a natural advantage in solving quasiperiodic problems with such Fourier coefficient distribution.
Compared with the case $v_0=2.5$, the concentrated area of the case $v_0=3$ is greatly elongated. 
Entire size of this concentrated area can reach $5082\times8192$.
Such a large computing area could be unaffordable for PM. 
However, by using the parallelogram index set $\mathcal{K}_{K,L}$ with a small $K$, IWFPM can still efficiently solve this case.
To demonstrate the superiority of our algorithm in handling the above two cases, we use both PM and IWFPM to solve this QSE.

\begin{table}[!hbpt]
    \centering
    \footnotesize{
    \caption{Condition numbers of $\mathbf{\tilde{H}}$ and  $\bm{M}\mathbf{\tilde{H}}$ for 1D QSE with potential \cref{eq:E1_equation} ($v_0=2.5$).}
    \label{tab:E1_condition}
    \begin{tabular}{|c|c|c|c|c|c|c|c|}
    \hline
    \multicolumn{2}{|c|}{} & \multicolumn{6}{c|}{Condition Number} 
    \\ \hline
    \multicolumn{2}{|c|}{DOF} & \multicolumn{3}{c|}{$\mathbf{\tilde{H}}$} & \multicolumn{3}{c|}{$\bm{M}\mathbf{\tilde{H}}$}
    \\ \hline
    \multicolumn{2}{|c|}{$L$} & 200 & 400 & 800 & 200 & 400 & 800\\
    \hline
    \multirow{2}*{$K=200$} & PM & 4.72e+05 & 9.02e+05 & 2.17e+06 & 2.24 & 2.24 & 2.24\\
    \cline{2-8}
    & IWFPM & 1.81e+05 & 1.81e+05 & 1.81e+05 & 2.24& 2.24&2.24\\
    \hline
    \multicolumn{2}{|c|}{$K$} & 50 & 100 & 200 & 50 & 100 & 200\\
    \hline
    \multirow{2}*{$L=1600$} & PM &4.87e+06& 5.35e+06 & 6.37e+06 &2.24& 2.24&2.24\\
    \cline{2-8}
    & IWFPM &1.15e+04&4.55e+04&1.81e+05 &2.24& 2.24&2.24\\
    \hline
    \end{tabular}
    }
\end{table}

\vspace{2mm}
\textbf{Case 1: $v_0=2.5$.}
First, in \Cref{tab:E1_condition}, we give a comparison of the condition numbers of $\mathbf{\tilde{H}}$ before and after preconditioning, to show the effectiveness of the preconditioner $\bm{M}$ defined by \cref{eq:precond}. 
The results show that whether using PM or IWFPM, the condition numbers of $\mathbf{\tilde{H}}$ generally exceed the magnitude of 1.0e+04. While after preconditioning, they are remarkably reduced from $>$1.0e+04 to 2.24. 
In fact, the preconditioner we designed has shown amazing condition number optimization effects for solving 1D, 2D, 3D QSEs and under different quantum states.

After using the powerful preconditioning, we next compare the algorithm accuracy of PM and IWFPM in terms of eigenvalue and eigenfunction errors. 
To present the tiny errors clearly, we set the calculated eigenvalue $E^*_0$ and eigenfunction $\bm{u}^*_0$ of a large-scale system as the numerical exact solution for comparison. 
Here, $E^*_0$ and $\bm{u}^*_0$ are calculated by PM with $K=640,~L=1024$. \cref{tab:E1_error1} records the errors of eigenvalue $E_0$ and eigenfunction $\bm{u}_0$, respectively. 
The data shows that both PM and IWFPM methods can achieve high accuracy in calculating this example. While, from this comparison, it can be seen that under the same $L$, PM will need a more lager $K$ to achieve the same convergence accuracy compared with IWFPM. 
For instance, when achieving $E_v=2.84$e-14 and $E_f=1.49$e-05, IWFPM requires $K=5$, while PM requires $K=0.7L=42$. 

\begin{table}[!hbpt]
    \centering
    \footnotesize{
        \caption{Errors of PM and IWFPM when solving 1D QSE with potential \cref{eq:E1_equation} ($v_0=2.5$).}
    \label{tab:E1_error1}
    \begin{tabular}{|c|c|c|c|c|c|}
    \hline
    \multicolumn{2}{|c|}{} &$L$ & 20 & 40&60\\
	\hline
    \multirow{4}*{$E_v$}&\multirow{3}*{PM}&$K=0.3L$&1.45e-05&7.35e-08&9.42e-10\\
    \cline{3-6}
    &&$K=0.5L$&4.60e-07&2.79e-10&4.13e-13\\
    \cline{3-6}
    &&$K=0.7L$&5.64e-08&2.69e-11&2.82e-14\\
    \cline{2-6}
    &IWFPM&$K=5$&5.64e-08&2.11e-11&2.84e-14\\
    \hline
    \multirow{4}*{$E_f$}&\multirow{3}*{PM}&$K=0.3L$&1.11e-01&1.91e-02&2.59e-03\\
    \cline{3-6}
    &&$K=0.5L$&3.00e-02&2.05e-03&5.13e-05\\
    \cline{3-6}
    &&$K=0.7L$&1.82e-02&4.03e-04&1.49e-05\\
    \cline{2-6}
    &IWFPM&$K=5$&1.82e-02&3.48e-04&1.49e-05\\
    \hline
    \end{tabular}
    }
\end{table}

\begin{figure*}[!hbpt]
    \centering
    \includegraphics[width=0.5\textwidth]{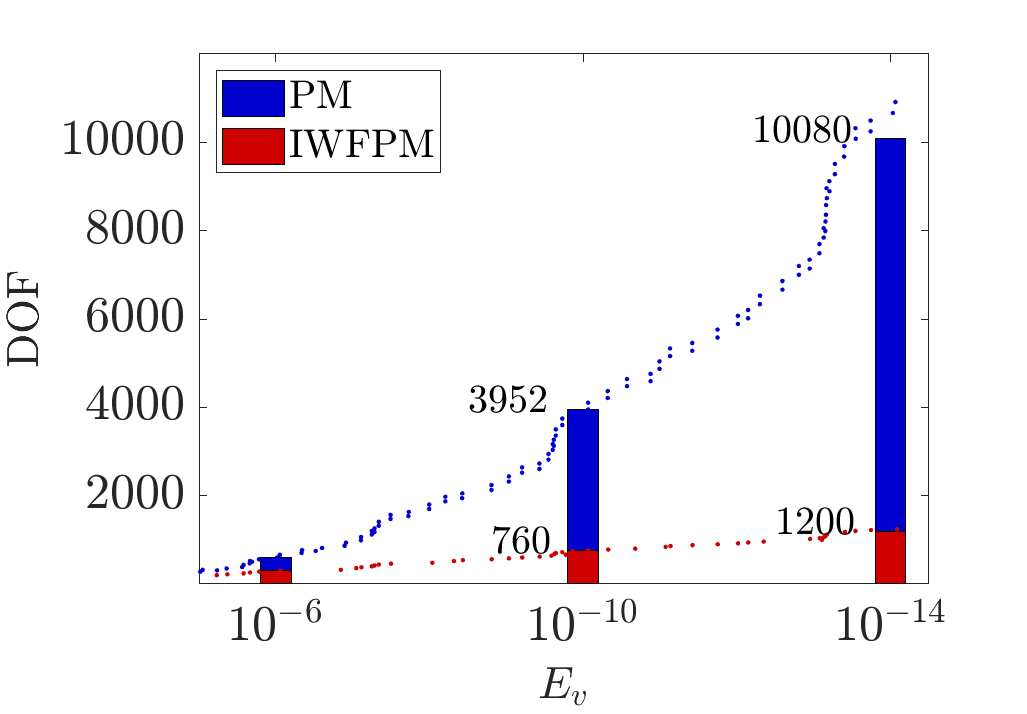}
    \caption{Required DOFs of PM and IWFPM when they achieve the same accurate $E_v$ for solving 1D QSE with potential \cref{eq:E1_equation} ($v_0=2.5$).} 
    \label{fig:E1_err_dof}
\end{figure*}

Note that the size of $K$ and $L$ affects the DOFs of the two algorithms, and fundamentally affects their calculation time. 
Hence, we next present \Cref{fig:E1_err_dof} to show the DOFs required by PM and IWFPM when they achieve the same errors of $E_0$. 
It is obvious that as the required accuracy increases, the DOF required for PM rises sharply compared with IWFPM.
When $E_v\approx2.8$e-14, the DOF of IWFPM is 1200, while the DOF of PM has exceeded 10000. 
We further compare the computational costs of PM and IWFPM. 
\Cref{tab:E1_CPU1} records the ITs and CPU times required by PM and IWFPM when $E_v\approx2.8$e-14.   
It is clear that although the ITs required by the two methods to achieve the same convergence accuracy are basically the same, the CPU time consumed differs by 10.75 times. 
In conclusion, the above series of results all demonstrate that IWFPM can greatly improve the calculation efficiency. 

\begin{table}[!hbpt]
    \centering
    \footnotesize{
    \caption{ITs and CPU times required by PM and IWFPM when solving the 1D QSE with potential \cref{eq:E1_equation} ($v_0=2.5$).}
    \label{tab:E1_CPU1}
    \begin{tabular}{|c|c|c|c|}
    \hline
    & PM & IWFPM & PM/IWFPM\\
    \hline
    $E_v$ &2.84e-14& 2.82e-14  & 1.01\\
    \hline
    DOF & 10080 & 1200 & 8.40\\
    \hline
    IT & 204 & 195 & 1.05\\
    \hline
    CPU(s) & 0.86 & 0.08 & 10.75\\
    \hline
    \end{tabular}
    }
\end{table}

\begin{table}[!hbpt]
    \centering
    \footnotesize{
        \caption{Results of PM and IWFPM for solving the 1D QSE with potential \cref{eq:E1_equation} ($v_0=3$).}
    \label{tab:E1_error_IWFPM}
    \begin{tabular}{|c|c|c|c|c|c|c|c|}
    \hline
    & \multicolumn{4}{c|}{IWFPM} & PM & PM/IWFPM\\
    \hline
    $E_v$ &
    9.78e-06 & 2.45e-06 & 6.03e-07 & \textbf{1.42e-07} & \textbf{1.44e-07} & \textbf{1.01}\\
    \hline
    $E_f$ & 
    8.43e-03 & 3.65e-03 & 4.16e-04 & \textbf{3.81e-04} &
    \textbf{3.57e-04} & \textbf{0.94}\\
    \hline
    DOF & 4096 & 8192 & 16384 & \textbf{32768} & \textbf{2621440} & \textbf{80}\\
    \hline
    IT & 1147 & 2140 & 4092 & \textbf{6871} & \textbf{7656} & \textbf{1.11}\\
    \hline
    CPU(s) & 0.57 & 1.72 & 6.08 & \textbf{22.87} & \textbf{3057.75} & \textbf{133.70}\\
    \hline
    \end{tabular}
    }
\end{table}

\vspace{2mm}
\textbf{Case 2: $v_0=3$.}
\Cref{tab:E1_error_IWFPM} records the results of PM and IWFPM for solving the eigenvalue $E_0$ and eigenfunction $\bm{u}_0$. 
Here, the numerical exact solutions $E^*_0$ and $\bm{u}^*_0$ are obtained by IWFPM with $K=10,~L=4096$. 
The reason why we not use PM to obtain a numerical reference solution is that the excessively large DOF in this case makes PM unaffordable. 
As we can see, IWFPM can achieve the same accuracy as PM with much fewer DOF, and less computational cost. 
For example, when the error $E_0\approx1.4$e-07, the DOF of IWFPM is $32768$, while PM is  $2621440$. 
It means that IWFPM can speed up the CPU time of PM by 133.70 times.

\begin{example}
    Consider 2D QSE \cref{eq:QSE} with potential
    \begin{equation}\label{eq:E2_equation}
        v(\bm{x})=4-[\cos(\beta x_1)+2\cos(\beta x_2)+\cos(\beta x_1\cos\theta+\beta x_2\sin\theta)],
    \end{equation}
    where $\bm{x}=(x_1,x_2)^T,~\beta\in\mathbb{R},~\theta\in(0,2\pi)$. 
\end{example}

The projection matrix of $v(\bm{x})$ is 
\begin{equation*}
    \bm{P}=\beta\begin{pmatrix}1&0&\cos\theta\\0&1&\sin\theta\end{pmatrix},
\end{equation*}
and the corresponding parent function is $$V(\bm{y})=4-(\cos y_1+2\cos y_2+\cos y_3),\qquad\bm{y}=(y_1,y_2,y_3)^T.$$
The parallelogram index set $\mathcal{K}_{K,L}$ defined by \cref{eq:non-tensor} is
\begin{footnotesize}
    $$
    \mathcal{K}_{K,L}=\left\{\bm{k}=(k_1,k_2,k_3)^T\in\mathbb{Z}^3: k_1+ k_3\cos\theta\in[-K,K), k_2+ k_3\sin\theta \in[-K,K), k_3\in[-L,L) \right\}.
    $$
\end{footnotesize}

\begin{figure*}[!hbpt]
    \centering
    \subfigure[I:  $\beta=0.8\pi,\theta=0.2\pi$]{\label{fig:E2_wave1}\includegraphics[width=0.3\textwidth]{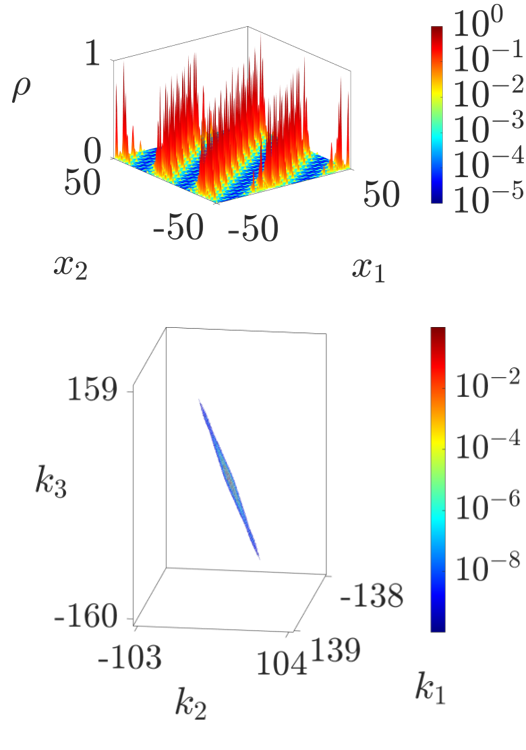}}
    \subfigure[II: $\beta=0.5\pi,\theta=0.25\pi$]{\label{fig:E2_wave2}\includegraphics[width=0.3\textwidth]{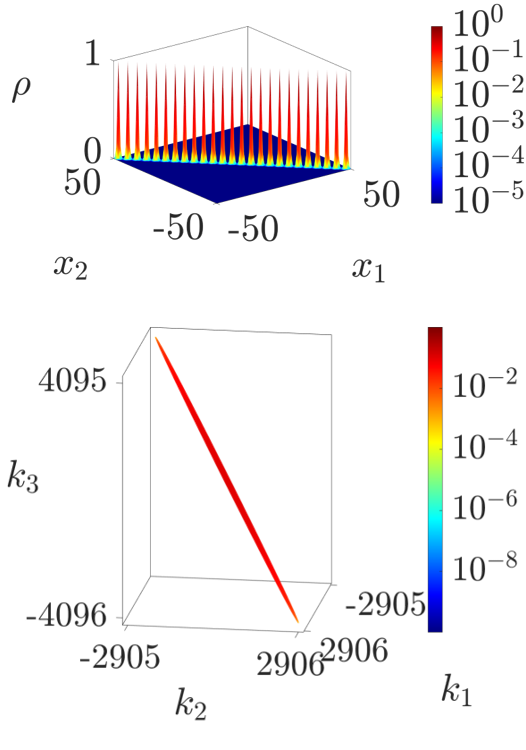}}
    \subfigure[III: $\beta=0.5\pi,\theta=0.2\pi$]{\label{fig:E2_wave3}\includegraphics[width=0.3\textwidth]{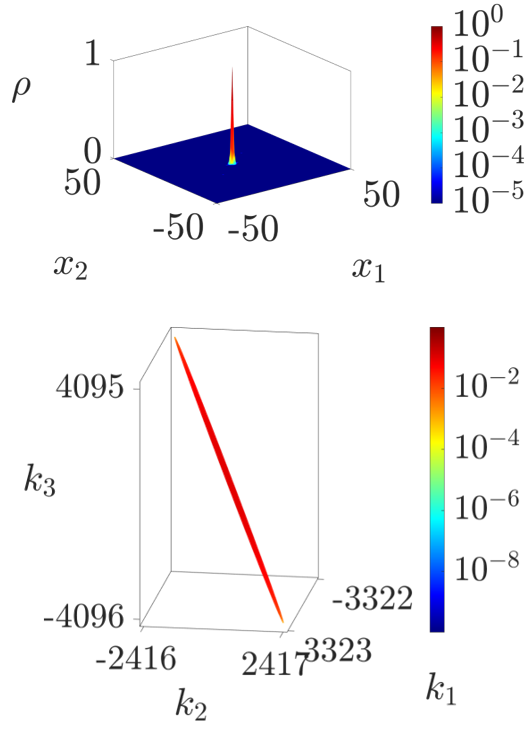}}
    \caption{Results of solving 2D QSE with potential \cref{eq:E2_equation} by IWFPM. The top row: probability density function $\bm{\rho}$; The bettom row:  Fourier coefficients $\tilde{U}_{\bm{k}}~(\tilde{U}_{\bm{k}}\geq10^{-8})$.}
    \label{fig:E2_wave}
\end{figure*}

\Cref{fig:E2_wave} shows the probability density function $\bm{\rho}$ under the potential \cref{eq:E2_equation} with different parameters $\beta$ and $\theta$. 
Among them, \Cref{fig:E2_wave1} shows a 2D extended state and \Cref{fig:E2_wave3} shows a 2D localized state. 
The phase transition between the two states arising from $\beta$ plays a dominant role in interfering in the degree of localization of the wave function.
While some special values of $\theta$ can bring periodicity to the wave function. 
As shown in \Cref{fig:E2_wave2}, when $\beta=0.5\pi,~\theta=0.25\pi$, the $\bm{\rho}$ exhibits extended state along the line $x_1+x_2=0$ and localized state in the orthogonal direction. 


Next, we consider the concentrated area of the Fourier coefficient $\tilde{U}_{\bm{k}}$ ($\gg$1.0e-8) under the potential \cref{eq:E2_equation} with different parameters, as shown in \cref{fig:E2_wave}. 
Similar to \Cref{exm:1}, localized states require more DOFs than that of extended states.
In the following, we use IWFPM to solve the above mentioned three quantum states and compare the results with PM. 

\vspace{2mm}
\textbf{Case 1: I: $\beta=0.8\pi,~\theta=0.2\pi$.} 
\cref{tab:E2_error1} records the errors of eigenvalue $E_0$ and eigenfunction $\bm{u}_0$, respectively. 
Here, the numerical exact solution $E^*_0$ and $\bm{u}^*_0$ are calculated by IWFPM when $K=10,~L=160$. 
Compared with the PM, IWFPM exhibits higher-order convergence under a small scale of $K$. 
Specifically, when $K=6,~L=20$, the error $E_v$ of IWFPM reaches 4.53e-08, while that of PM does not reach 7.51e-05. To achieve the same magnitude of error, PM needs $K=16,~L=20$. 

\begin{table}[!hbpt]
    \centering
    \footnotesize{
        \caption{Errors of PM and IWFPM when solving 2D QSE with potential \cref{eq:E2_equation} (I:  $\beta=0.8\pi,~\theta=0.2\pi$).}
    \label{tab:E2_error1}
    \begin{tabular}{|c|c|c|c|c|c|c|}
    \hline
    \multicolumn{2}{|c|}{} &$L$ & 20 & 30&40&50\\
	\hline
    \multirow{4}*{$E_v$} & \multirow{3}*{PM}
    &$K=0.4L$&7.51e-05&2.55e-06&6.21e-08&1.83e-09\\
    \cline{3-7}
    & &$K=0.6L$&2.55e-06&1.48e-08&8.60e-11&1.65e-12\\
    \cline{3-7}
    & &$K=0.8L$&7.35e-08&9.88e-11&3.85e-13&1.11e-15\\
    \cline{2-7}
    & IWFPM&$K=6$&4.53e-08&3.09e-11&1.84e-13&1.11e-15\\
    \hline
    \multirow{4}*{$E_f$} & \multirow{3}*{PM} 
    &$K=0.4L$&1.55e-01&2.58e-02&3.48e-03&7.99e-04\\
    \cline{3-7}
    &&$K=0.6L$&2.58e-02&1.94e-03&2.22e-04&3.46e-05\\
    \cline{3-7}
    &&$K=0.8L$&3.81e-03&2.35e-04&1.43e-05&1.18e-06\\
    \cline{2-7}
    &IWFPM&$K=6$&3.12e-03&1.67e-04&1.00e-05&1.09e-06\\
    \hline
    \end{tabular}
    }
\end{table}
\begin{figure*}[!hbpt]
    \centering
    \includegraphics[width=0.5\textwidth]{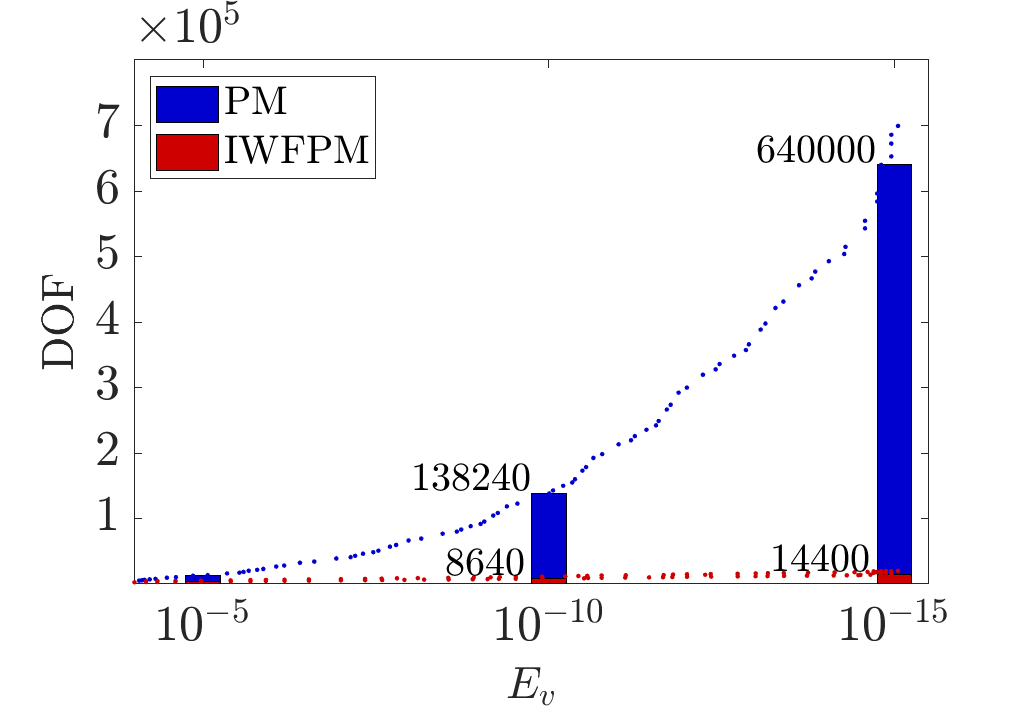}
    \caption{Required DOFs by PM and IWFPM when they arrive in the same accurate $E_v$ for solving 2D QSE with potential \cref{eq:E2_equation} (I:  $\beta=0.8\pi,~\theta=0.2\pi$).} 
    \label{fig:E2_err_dof}
\end{figure*}

\begin{table}[!hbpt]
    \centering
    \footnotesize{
    \caption{Comparison of ITs and CPU times spent by PM and IWFPM, when solving the 2D QSE with potential \cref{eq:E2_equation} (I:  $\beta=0.8\pi,~\theta=0.2\pi$) with $E_v\approx1.1$e-15.}
    \label{tab:E2_CPU1}
    \begin{tabular}{|c|c|c|c|}
    \hline
    & PM & IWFPM & PM/IWFPM\\
    \hline
    $E_v$ & 1.11e-15 & 1.11e-15 & 1.00\\
    \hline
    DOF & 640000 & 14400 & 44.44\\
    \hline
    IT & 227 & 229 & 0.99\\
    \hline
    CPU(s) & 25.27 & 0.39 & 64.79\\
    \hline
    \end{tabular}
    }
\end{table}

We further compare the ITs and CPU times of PM and IWFPM when achieving the same accuracy. 
First, \Cref{fig:E2_err_dof} illustrates the DOFs required by PM and IWFPM when they achieve the same errors of $E_0$. Once again, the results underscore that, in comparison to IWFPM, PM exhibits a notable drawback in terms of computational storage. 
Then, we consider the case when the error is about 1.1e-15. The data in \Cref{tab:E2_CPU1} shows that although PM and IWFPM have almost the same IT, IWFPM can reduce the DOF of PM by 44 times, which ultimately results in IWFPM taking over 64 times CPU time less than PM.

\vspace{2mm}
\textbf{Case 2: II: $\beta=0.5\pi,~\theta=0.25\pi$ and III: $\beta=0.5\pi,~\theta=0.2\pi$.}
\Cref{tab:E2_2_error_IWFPM} and \Cref{tab:E2_3_error_IWFPM} present the errors, ITs, and CPU times of IWFPM and PM for solving the eigenvalue $E_0$ and eigenfunction $u_0$. 
Here, the numerical exact solution $E^*_0$ and $u^*_0$ are calculated by IWFPM when $K=10,~L=4096$.
The data once again shows the high efficiency and accuracy of IWFPM, regardless of whether the eigenstate is periodic localized or localized. 
When the two methods both achieve the error $E_v\approx1.5$e-05, it is apparent that the DOF required by PM is 156 times larger than that needed by IWFPM. 
As a result, the CPU time consumed by PM is 168 times larger than that by IWFPM. 
Note that, if higher accuracy is achieved, the advantages of IWFPM over PM will be even more exaggerated.

\begin{table}[!hbpt]
    \centering
    \footnotesize{
        \caption{Results of PM and IWFPM for solving the 2D QSE with potential \cref{eq:E2_equation} (II: $\beta=0.5\pi,~\theta=0.25\pi$).}
    \label{tab:E2_2_error_IWFPM}
    \begin{tabular}{|c|c|c|c|c|c|c|c|}
    \hline
    & \multicolumn{4}{c|}{IWFPM} & PM & PM/IWFPM\\
    \hline
    $E_v$ & \textbf{1.51e-05} & 3.76e-06 & 9.30e-07 & 2.22e-07 & \textbf{1.51e-05} & \textbf{1.00}\\
    \hline
    $E_f$ & \textbf{4.34e-03} & 1.08e-03 & 2.88e-04 & 5.64e-05 & \textbf{4.36e-03} & \textbf{1.00}\\
    \hline
    DOF & \textbf{65536} & 131072 & 262144 & 524288 & \textbf{10240000} & \textbf{156.25}\\
    \hline
    IT & \textbf{1873} & 3538 & 6726 & 13133 & \textbf{1653} & \textbf{0.88}\\
    \hline
    CPU(s) & \textbf{16.98} & 71.56 & 264.67 & 1009.56 & \textbf{3050.37} & \textbf{179.64}\\
    \hline
    \end{tabular}
    }
\end{table}
\begin{table}[!hbpt]
    \centering
    \footnotesize{
        \caption{Results of PM and IWFPM for solving the 2D QSE with potential \cref{eq:E2_equation} (III: $\beta=0.5\pi,~\theta=0.2\pi$).}
    \label{tab:E2_3_error_IWFPM}
    \begin{tabular}{|c|c|c|c|c|c|c|c|}
    \hline
    & \multicolumn{4}{c|}{IWFPM} & PM & PM/IWFPM\\
    \hline
    $E_v$ & \textbf{1.43e-05} & 3.58e-06 & 8.85e-07 & 2.11e-07 & \textbf{1.54e-05} & \textbf{1.08}\\
    \hline
    $E_f$ & \textbf{9.59e-01} & 9.01e-01 & 6.94e-01 & 1.75e-01 & \textbf{9.61e-01} & \textbf{1.00}\\
    \hline
    DOF &  \textbf{65536} & 131072 & 262144 & 524288 & \textbf{10240000} & \textbf{156.25}\\
    \hline
    IT & \textbf{1911} & 3612 & 6870 & 13426 & \textbf{1653} & \textbf{0.87}\\
    \hline
    CPU(s) & \textbf{18.53} & 77.50 & 274.81 & 1025.83 & \textbf{3111.78} & \textbf{167.93}\\
    \hline
    \end{tabular}
    }
\end{table}

\begin{example}
    Consider 3D QSE \cref{eq:QSE} with potential
    \begin{equation}\label{eq:E4_equation}
    \begin{aligned}
        v(\bm{x})=6-[&\cos(\beta x_1)+\cos(\beta x_2)+\cos(\beta x_3)+\cos(\beta x_1\cos\theta+\beta x_2\sin\theta))\\
        &+\cos(-\beta x_1\sin\theta+\beta x_2\cos\theta)+\cos(\beta \alpha x_3)]
    \end{aligned}
    \end{equation}
    where $\bm{x}=(x_1,x_2,x_3)^T,~\beta\in\mathbb{R},~\theta=0.2\pi,~\alpha=(\sqrt{5}-1)/2$. 
\end{example}
The projection matrix of $v(\bm{x})$ is 
\begin{equation*}
\bm{P}=\beta\begin{pmatrix}
1&0&0&\cos\theta&-\sin\theta&0\\
0&1&0&\sin\theta&\cos\theta&0\\
0 & 0 & 1 & 0 & 0 & \alpha
\end{pmatrix},
\end{equation*}
and the corresponding parent function is $$\mathcal{V}(\bm{y})=6-\cos \bm{y},\qquad\bm{y}=(y_1,\cdots,y_6)^T.$$
The spectral point set $\mathcal{K}_{K,L}$ defined by \cref{eq:non-tensor} is 
\begin{equation*}
\begin{aligned}
\mathcal{K}_{K,L}=\{\bm{k}= (k_1,\cdots,k_6)^T\in \mathbb{Z}^6: k_1+k_4 \cos\theta - k_5\sin\theta  \in[-K,K), & \\ 
k_2+ k_4 \sin\theta + k_5\cos\theta \in[-K,K),~k_3+ k_6 \alpha \in[-K,K),&\\ 
k_4\in[-L,L),~k_5\in[-L,L),~k_6\in[-L,L)\}.&  
\end{aligned}
\end{equation*}

\begin{figure*}[!hbpt]
    \centering
    \subfigure[$\rho$]{\label{fig:E4_wave1}\includegraphics[width=0.39\textwidth]{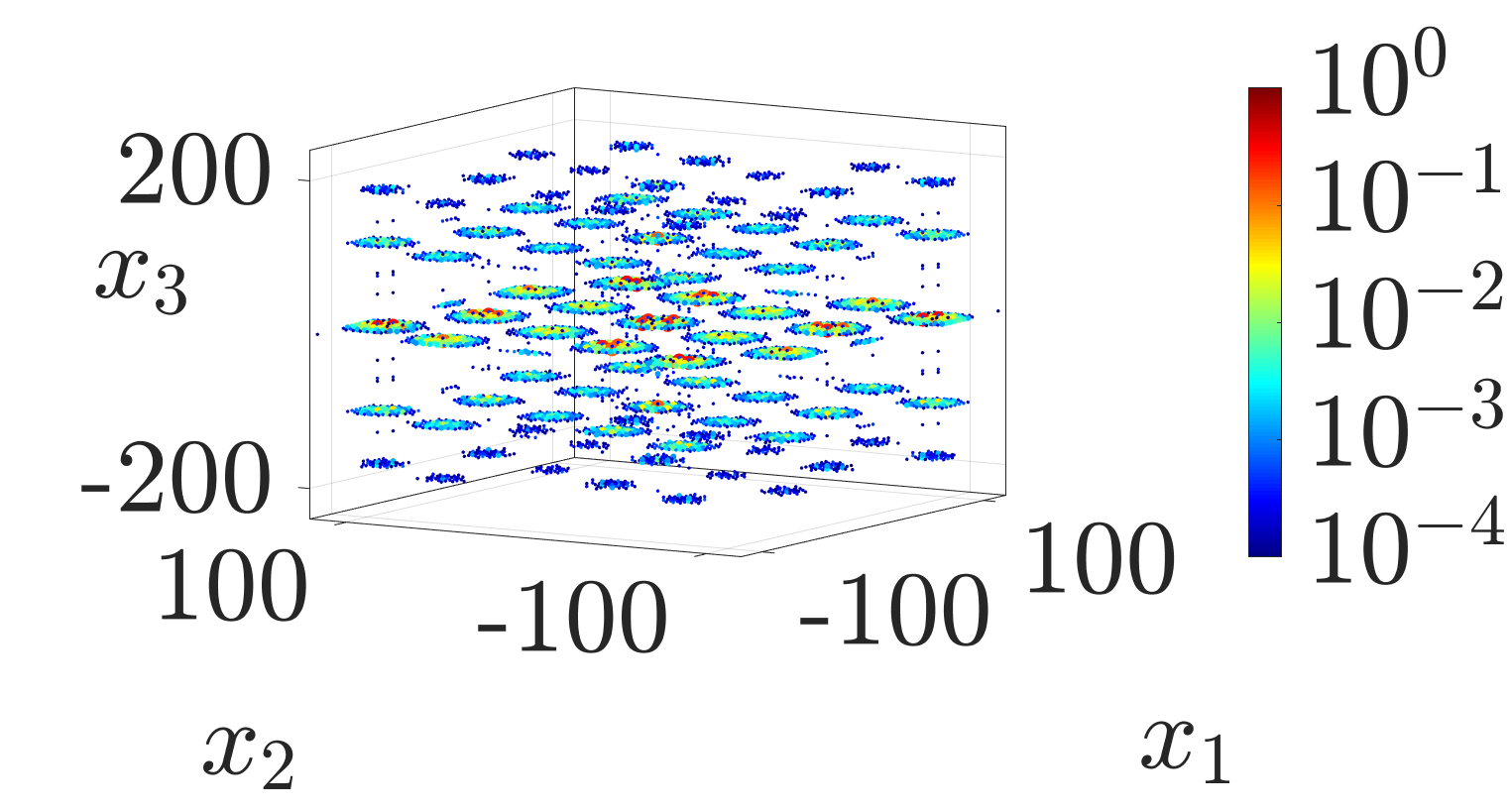}}
    \subfigure[$\tilde{U}_{\bm{k}}$($k_4,k_5,k_6=0$)]{\label{fig:E4_spectral1_s0}\includegraphics[width=0.29\textwidth]{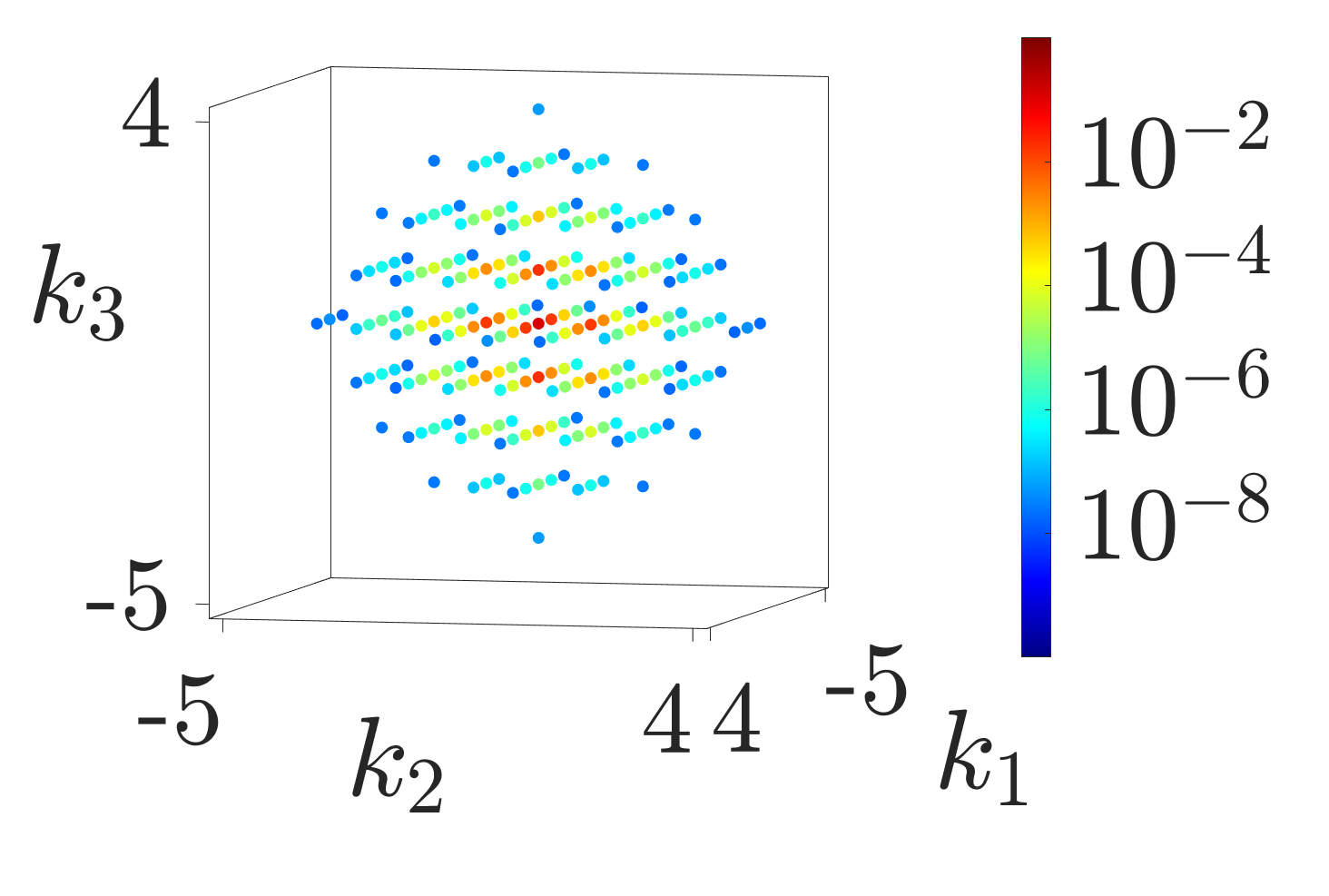}}
    \subfigure[$\tilde{U}_{\bm{k}}$($k_4,k_5,k_6=12$)]{\label{fig:E4_spectral1_s12}\includegraphics[width=0.29\textwidth]{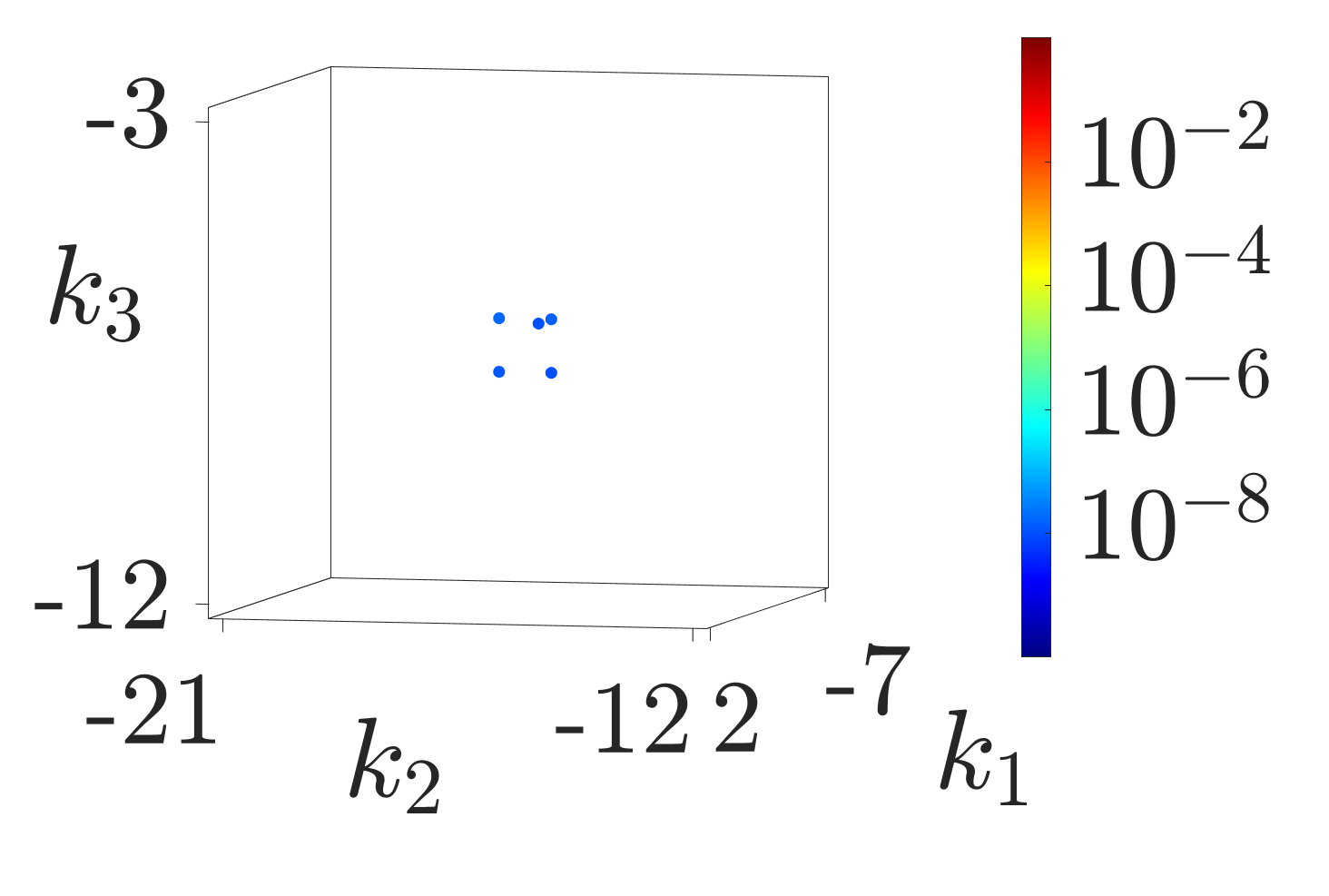}}
    \caption{Results of solving 3D QSE with potential \cref{eq:E4_equation} by IWFPM ($\beta=\pi$). (a) The probability density function $\bm{\rho}$; (b)(c) Slices of the Fourier coefficients $\tilde{U}_{\bm{k}}~(\tilde{U}_{\bm{k}}\geq10^{-8})$.}
    \label{fig:E4_1}
 \end{figure*}
 
\begin{figure*}[!hbpt]
    \centering
    \subfigure[$\rho$]{\label{fig:E4_wave2}\includegraphics[width=0.39\textwidth]{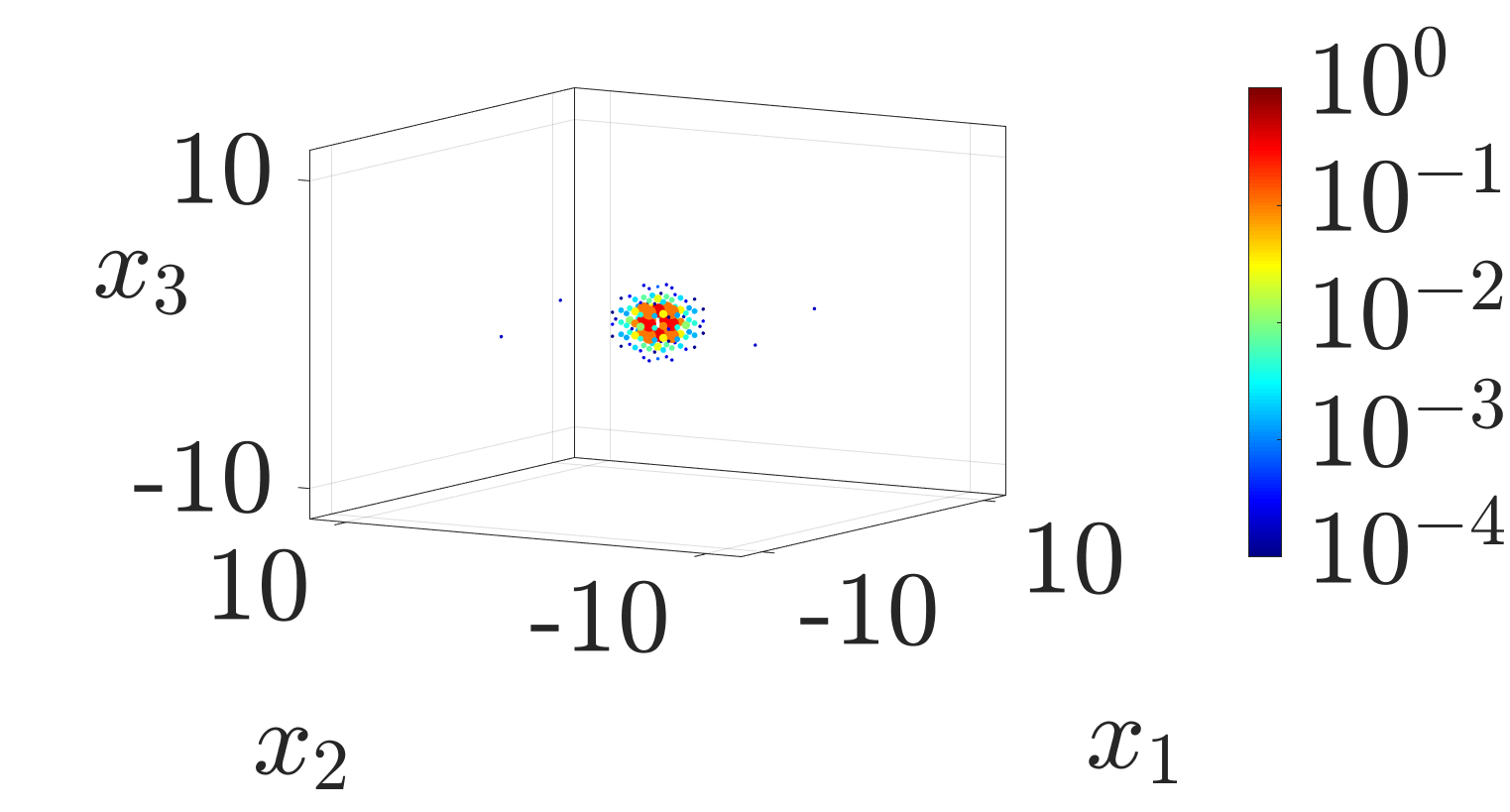}}
    \subfigure[$\tilde{U}_{\bm{k}}$($k_4,k_5,k_6=0$)]{\label{fig:E4_spectral2_s0}\includegraphics[width=0.29\textwidth]{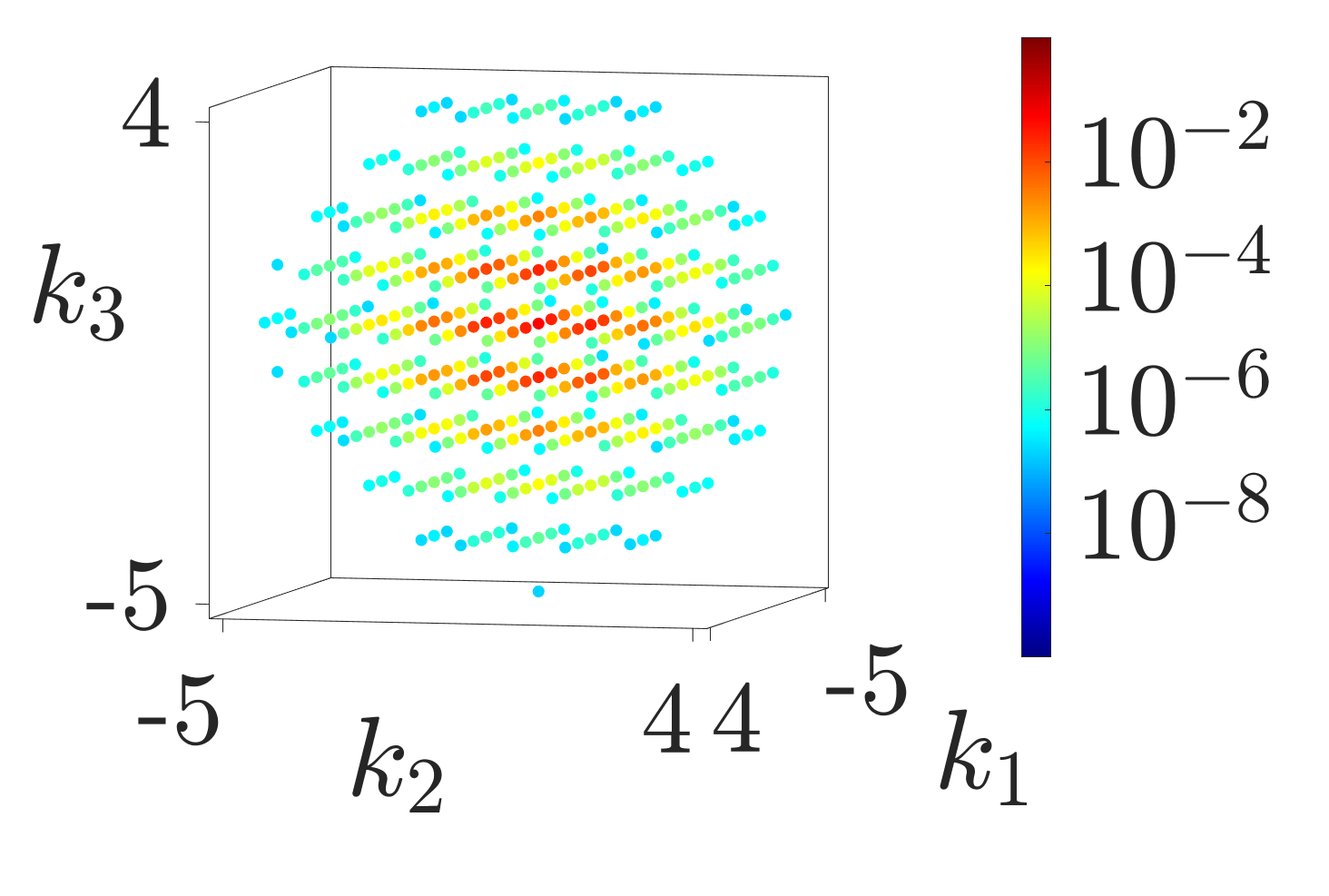}}
    \subfigure[$\tilde{U}_{\bm{k}}$($k_4,k_5,k_6=32$)]{\label{fig:E4_spectral2_s32}\includegraphics[width=0.29\textwidth]{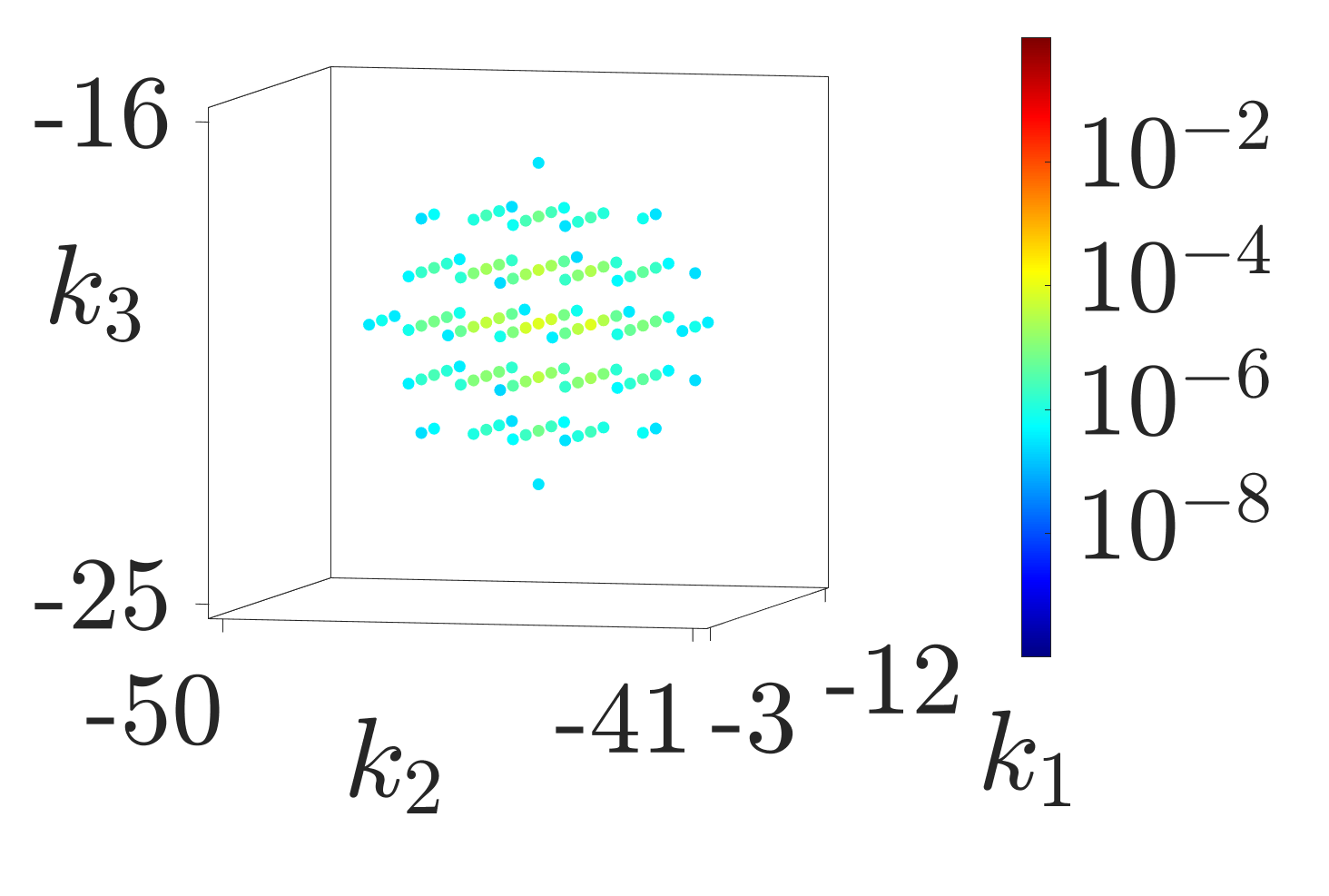}}
    \caption{Results of solving 3D QSE with potential \cref{eq:E4_equation} by IWFPM ($\beta=0.5\pi$). (a) The probability density function $\bm{\rho}$; (b)(c) Slices of the Fourier coefficients $\tilde{U}_{\bm{k}}~(\tilde{U}_{\bm{k}}\geq10^{-8})$.}
     \label{fig:E4_2}
 \end{figure*}

\Cref{fig:E4_wave1} and \Cref{fig:E4_wave2} show the probability density function $\rho$ under the quasiperiodic potential \cref{eq:E4_equation} with different parameters $\beta$. 
When $\beta=\pi$, the wave function diffuses throughout the entire 3D space, which seems to be a 3D extended state. 
Conversely, when $\beta=0.5\pi$, the wave function becomes localized eigenstate. 

To show the concentrated distribution of Fourier coefficients within a 3D parallelogram area, we depict the slices of $\tilde{U}_{\bm{k}}$ under the potential \cref{eq:E4_equation} with different parameters in \Cref{fig:E4_1} and \Cref{fig:E4_2}. 
Clearly, the $L$ required for the case I is less than 12, while the $L$ required for the case II is more than 32. 
To calculate the case II in the 6D space will bring a huge storage difficulty for PM (DOF exceeding $10^8$). 
However, IWFPM enables the calculation of this problem by significantly reducing DOFs. 
\Cref{tab:E4_error5} records the results of IWFPM for solving the eigenvalue $E_0$ and eigenfunction $\bm{u}_0$. 
Here, the numerical exact solutions $E^*_0$ and $\bm{u}^*_0$ of $E_0$ and $\bm{u}_0$ are obtained by IWFPM with $K=5,~L=40$. 
The data show that IWFPM can achieve convergence when calculating 3D QSE with \cref{eq:E4_equation}, regardless of whether the solution belongs to extended state or localized state.

\begin{table}[!hbpt]
    \centering
    \footnotesize{
    \caption{Results of IWFPM for solving the 3D QSE with potential \cref{eq:E4_equation} (I: $\beta=\pi$; II: $\beta=0.5\pi$).}
    \label{tab:E4_error5}
    \begin{tabular}{|c|c|c|c|c|c|c|c|c|c|c|}
    \hline
    $K=5$ & \multicolumn{2}{c|}{$E_v$} & \multicolumn{2}{c|}{$E_f$}\\
    \hline
    $L$ & I & II & I & II\\
    \hline
    10 & 6.82e-04 & 4.08e-03 & 2.88e-04 & 3.55e-02\\
    \hline
    20 & 1.58e-04 & 9.91e-04 & 1.94e-04 & 4.84e-03\\
    \hline
    30 & 4.34e-05 & 2.71e-04 & 9.32e-05 & 2.21e-03\\
    \hline
    \end{tabular}
    }
\end{table}

\section{Conclusion and outlook}\label{sec:concl}

In this paper, a new algorithm IWFPM is proposed based on PM for quasiperiodic systems with concentrated spectral point distribution. 
It filters out dominant spectral points by defining an irrational window and uses a corresponding index-shift transform to make the FFT available. 
The convergence analysis and computational cost of IWFPM are also given. 
We apply IWFPM to 1D, 2D, and 3D QSEs to demonstrate its accuracy and efficiency. 
An efficient diagonal preconditioner is also designed for the discrete QSEs to significantly reduce condition number. 
For both extended and localized quantum states, IWFPM exhibits a significant computational advantage over PM. 
More importantly, by using IWFPM, the existence of Anderson localization in 2D and 3D QSEs is numerically verified.

The proposed method opens avenues for further research in several directions. Firstly, given the widespread existence of the observed spectral point distribution feature, it is imperative to apply the IWFPM to a broader range of quasiperiodic systems to demonstrate its applicability. Secondly, the convergence analysis of quasiperiodic eigenproblems differs from conventional numerical analysis, as the Hilbert space for describing potential energy functions and eigenfunctions is entirely different. New fundamental theories of numerical analysis are needed to describe this problem, and we plan to explore this in our future work.
Thirdly, we will extend the application of our method to more quasiperiodic systems, aiming to discover exotic phenomena and even unveil new physical principles.

\appendix
\section{Proof of \Cref{lm:norm relation}}
\label{sec:proof_norm_relation}
\begin{proof}
We only prove the equivalence $|u|_{\alpha,\beta}\simeq |u|_{\alpha}+|{U}|_{\mathcal{H}^\beta}$, since the other is similar. 
First, it is easy to obtain $|u|_{\alpha,\beta}\lesssim |u|_{\alpha}+|{U}|_{\mathcal{H}^\beta}$ since
\begin{equation*}
    \begin{aligned}
    |u|_{\alpha,\beta}^2&=\sum_{\bm{k}\in\mathbb{Z}^n}\big(\|\bm{k}_{\mathrm{I}}+\bm{Qk}_{\mathrm{I}}\|^{2\alpha}+\|\bm{k}_{\mathrm{II}}\|^{2\beta}\big)|\hat{{U}}_{\bm{k}}|^2\\
    &=\sum_{\bm{k}\in\mathbb{Z}^n}\big(\|\bm{P}_{\mathrm{I}}^{-1}\|^{2\alpha}\|\bm{Pk}\|^{2\alpha}+ \|\bm{k}\|^{2\beta}\big)|\hat{{U}}_{\bm{k}}|^2\\
    &\lesssim\sum_{\bm{k}\in\mathbb{Z}^n}\|\bm{Pk}\|^{2\alpha}|\hat{{U}}_{\bm{k}}|^2+\sum_{\bm{k}\in\mathbb{Z}^n}\|\bm{k}\|^{2\beta}|\hat{{U}}_{\bm{k}}|^2.
    \end{aligned}
\end{equation*}

Next, we prove the converse $|u|_{\alpha}+|{U}|_{\mathcal{H}^\beta}\lesssim |u|_{\alpha,\beta}$. It can be divided into two parts: $|u|_{\alpha}\lesssim|u|_{\alpha,\beta}$ and $|{U}|_{\mathcal{H}^\beta}\lesssim|u|_{\alpha,\beta}$. The former can be obtained naturally because 
\begin{equation*}
    \begin{aligned}
    |u|_{\alpha}^2=\sum_{\bm{k}\in\mathbb{Z}^n}\|\bm{Pk}\|^{2\alpha}|\hat{{U}}_{\bm{k}}|^2
    =\sum_{\bm{k}\in\mathbb{Z}^n}\|\bm{P}_{\mathrm{I}}\|^{2\alpha}\|\bm{k}_{\mathrm{I}} +\bm{Qk}_{\mathrm{II}}\|^{2\alpha}|\hat{{U}}_{\bm{k}}|^2.
    \end{aligned}
\end{equation*}
For the latter
\begin{equation*}
    |{U}|^2_{\mathcal{H}^\beta}=\sum_{\bm{k}\in\mathbb{Z}^n}\|\bm{k}\|^{2\beta}|\hat{{U}}_{\bm{k}}|^2=\sum_{\bm{k}\in\mathbb{Z}^n}\big(\|\bm{k}_{\mathrm{I}}\|^{2\beta}+\|\bm{k}_{\mathrm{II}}\|^{2\beta}\big)|\hat{{U}}_{\bm{k}}|^2, 
\end{equation*}
then we only need to focus on $\|\bm{k}_{\mathrm{I}}\|^{2\beta}$. 
We discuss it from three cases.
The first case is $\bm{k}_{\mathrm{II}}=0$, and it is obvious
$$\|\bm{k}_{\mathrm{I}}\|^{2\beta}=\|\bm{k}_{\mathrm{I}} +\bm{Qk}_{\mathrm{II}}\|^{2\beta}\leq\|\bm{k}_{\mathrm{I}} +\bm{Qk}_{\mathrm{II}}\|^{2\alpha}.$$
The second case is $\bm{k}_{\mathrm{II}}\neq0$ and $\|\bm{k}_{\mathrm{I}} +\bm{Qk}_{\mathrm{II}}\|<1$, and it yields
$$\|\bm{k}_{\mathrm{I}}\|\leq\|\bm{k}_{\mathrm{I}}+\bm{Qk}_{\mathrm{II}}\|+\|\bm{Qk}_{\mathrm{II}}\|< 1+\|\bm{Qk}_{\mathrm{II}}\|\leq \big(1+\|\bm{Q}\|\big)\|\bm{k}_{\mathrm{II}}\|.$$
The third case is $\bm{k}_{\mathrm{II}}\neq0$ and $\|\bm{k}_{\mathrm{I}} +\bm{Qk}_{\mathrm{II}}\|\geq1$, and we have
$$\|\bm{k}_{\mathrm{I}}\|\leq\|\bm{k}_{\mathrm{I}} +\bm{Qk}_{\mathrm{II}}\|+ \|\bm{Q}\|\|\bm{k}_{\mathrm{II}}\|,$$
which yields
$$\begin{aligned}
\|\bm{k}_{\mathrm{I}}\|^{2\beta}\leq C_1\left(\|\bm{k}_{\mathrm{I}} +\bm{Qk}_{\mathrm{II}}\|^{2\beta}+ \|\bm{k}_{\mathrm{II}}\|^{2\beta}\right)
\leq C_1\left(\|\bm{k}_{\mathrm{I}} +\bm{Qk}_{\mathrm{II}}\|^{2\alpha}+ \|\bm{k}_{\mathrm{II}}\|^{2\beta}\right),
\end{aligned}$$
where $C_1$ is a constant independent of $\bm{k}$. The above inequalities can lead to
$$\|\bm{k}\|^{2\beta}\leq C_2\left(\|\bm{k}_{\mathrm{I}} +\bm{Qk}_{\mathrm{II}}\|^{2\alpha}+ \|\bm{k}_{\mathrm{II}}\|^{2\beta}\right),$$
and so $|{U}|_{\mathcal{H}^\beta}\lesssim|u|_{\alpha,\beta}$.
\end{proof}

\section{Proof of \Cref{th:proj}}
\label{sec:proof_proj}
\begin{proof}
First of all, it is easy to obtain the following inequality
$$\begin{aligned}
\left|u-\mathcal{P}_{K,L}u\right|_{\mu,\nu}^2&=\sum_{\bm{k}\in\mathbb{Z}^n\setminus\mathcal{K}_{K,L}}\big(\|\bm{k}_{\mathrm{I}} +\bm{Qk}_{\mathrm{II}}\|^{2\mu}+\|\bm{k}_{\mathrm{II}}\|^{2\nu}\big)|\hat{{U}}_{\bm{k}}|^2\\&\leq (a_1+a_2)+(a_3+a_4),
\end{aligned}$$
where
$$\begin{aligned}
a_1=&\sum_{\bm{k}_{\mathrm{II}}\in\mathbb{Z}^{n-d}}~\sum_{\|\bm{k}_{\mathrm{I}} +\bm{Qk}_{\mathrm{II}}\|_{\infty}\geq K} \|\bm{k}_{\mathrm{I}} +\bm{Qk}_{\mathrm{II}}\|^{2\mu}|\hat{{U}}_{\bm{k}}|^2,\\
a_2=&\sum_{\|\bm{k}_{\mathrm{II}}\|_{\infty}\geq L} ~ \sum_{\|\bm{k}_{\mathrm{I}} +\bm{Qk}_{\mathrm{II}}\|_{\infty}<K} \|\bm{k}_{\mathrm{I}} +\bm{Qk}_{\mathrm{II}}\|^{2\mu}|\hat{{U}}_{\bm{k}}|^2,\\
a_3=&\sum_{\bm{k}_{\mathrm{I}}\in\mathbb{Z}^d}~\sum_{\|\bm{k}_{\mathrm{II}}\|_{\infty}\geq L}\|\bm{k}_{\mathrm{II}}\|^{2\nu}|\hat{{U}}_{\bm{k}}|^2,\\
a_4=&\sum_{\|\bm{k}_{\mathrm{II}}\|_{\infty}< L}~ \sum_{\|\bm{k}_{\mathrm{I}} +\bm{Qk}_{\mathrm{II}}\|_{\infty}\geq K} \|\bm{k}_{\mathrm{II}}\|^{2\nu}|\hat{{U}}_{\bm{k}}|^2,\\
\end{aligned}$$
and $\|\bm{m}\|_{\infty}=\max_{1\leq j\leq n'}\{|m_j|\}$ for $\bm{m}\in\mathbb{R}^{n'}$. Using the equivalence of vector norms, we have
$$\begin{aligned}
a_1=&\sum_{\bm{k}_{\mathrm{II}}\in\mathbb{Z}^{n-d}}~\sum_{\|\bm{k}_{\mathrm{I}} +\bm{Qk}_{\mathrm{II}}\|_{\infty}\geq K}\|\bm{k}_{\mathrm{I}} +\bm{Qk}_{\mathrm{II}}\|^{2\mu-2\alpha} \|\bm{k}_{\mathrm{I}}+\bm{Qk}_{\mathrm{II}}\|^{2\alpha} |\hat{{U}}_{\bm{k}}|^2\\
\leq&~K^{2\mu-2\alpha}\sum_{\bm{k}_{\mathrm{II}}\in\mathbb{Z}^{n-d}}~\sum_{\|\bm{k}_{\mathrm{I}} +\bm{Qk}_{\mathrm{II}}\|_{\infty}\geq K} \|\bm{P}_{\mathrm{I}}^{-1}\|^{2\alpha}\|\bm{P}\bm{k}\|^{2\alpha} |\hat{{U}}_{\bm{k}}|^2 \\
\lesssim&~K^{2\mu-2\alpha}|u|_{\alpha}^2,\\
a_2\leq&~K^{2\mu}\sum_{\|\bm{k}_{\mathrm{II}}\|_{\infty}\geq L} ~ \sum_{\|\bm{k}_{\mathrm{I}} +\bm{Qk}_{\mathrm{II}}\|_{\infty}<K}
\|\bm{k}_{\mathrm{II}}\|^{-2\beta}\|\bm{k}_{\mathrm{II}}\|^{2\beta}
|\hat{{U}}_{\bm{k}}|^2\lesssim K^{2\mu}L^{-2\beta}|{U}|_{\mathcal{H}^\beta}^2,\\
a_3=&\sum_{\bm{k}_{\mathrm{I}}\in\mathbb{Z}^d}~\sum_{\|\bm{k}_{\mathrm{II}}\|_{\infty}\geq L} \|\bm{k}_{\mathrm{II}}\|^{2\nu-2\beta}\|\bm{k}_{\mathrm{II}}\|^{2\beta}|\hat{{U}}_{\bm{k}}|^2\leq L^{2\nu-2\beta}|{U}|_{\mathcal{H}^\beta}^2,\\
a_4\leq&~L^{2\nu}\sum_{\|\bm{k}_{\mathrm{II}}\|_{\infty}< L}~ \sum_{\|\bm{k}_{\mathrm{I}} +\bm{Qk}_{\mathrm{II}}\|_{\infty}\geq K}\|\bm{Pk}\|^{-2\alpha}\|\bm{Pk}\|^{2\alpha} |\hat{{U}}_{\bm{k}}|^2\\
=&~L^{2\nu}\|\bm{P}_{\mathrm{I}}\|^{-2\alpha}\sum_{\|\bm{k}_{\mathrm{II}}\|_{\infty}< L}~ \sum_{\|\bm{k}_{\mathrm{I}} +\bm{Qk}_{\mathrm{II}}\|_{\infty}\geq K}\|\bm{\bm{k}_{\mathrm{I}} +\bm{Qk}_{\mathrm{II}}}\|^{-2\alpha}\|\bm{Pk}\|^{2\alpha} |\hat{{U}}_{\bm{k}}|^2\\
\lesssim&~L^{2\nu}K^{-2\alpha}|u|_{\alpha}^2.
\end{aligned}$$
Combining the above four equations, we can obtain the estimate of $\left|u-\mathcal{P}_{K,L}u\right|_{\mu,\nu}$ in the original proposition. Following the similar proof, we can also obtain the estimate of $\left\|u-\mathcal{P}_{K,L}u\right\|_{\mu,\nu}$.
\end{proof}

\section{Proof of \Cref{th:int}}
\label{sec:proof_int}
\begin{proof}
From \cref{eq:new_inner_product}, it is easy to prove that the Fourier basis functions satisfy the following discrete orthogonality
\begin{equation}\label{eq:discrete orthogonality}
\begin{aligned}
\left<\varphi_{\bm{k}},\varphi_{\bm{k}'}\right>_{K,L}=\left\{\begin{array}{ll}
         1, & \bm{k}-\bm{k}' \in\mathcal{Z}_{K,L}, \\
         0, & \text{otherwise},
       \end{array}
\right.
\end{aligned}
\end{equation}
where
$$\mathcal{Z}_{K,L}:=\{\bm{m}=(\bm{m}_{\mathrm{I}},\bm{m}_{\mathrm{II}})^T\in\mathbb{Z}^n:\bm{m}_{\mathrm{I}}/(2K)\in\mathbb{Z}^d,~ \bm{m}_{\mathrm{II}}/(2L)\in\mathbb{Z}^{n-d}\}.$$ 
This discrete orthogonality leads to the following aliasing formula
\begin{equation}\label{eq:d1 aliasing}
\tilde{u}_{\bm{k}}=\Bigg<\sum_{\bm{l}\in\mathbb{Z}^n}\hat{{U}}_{\bm{l}}\varphi_{\bm{Pl}},\varphi_{\bm{Pk}}\Bigg>_{K,L} =\sum_{\bm{m}\in\mathbb{Z}^n}\hat{{U}}_{(\bm{k}_{\mathrm{I}}+2K\bm{m}_{\mathrm{I}},\bm{k}_{\mathrm{II}}+2L\bm{m}_{\mathrm{II}})},
\end{equation}
for any $\bm{k}\in\mathcal{K}_{K,L}$.
Thus,
$$\mathcal{P}_{K,L}u(\bm{y})-\mathcal{I}_{K,L}u(\bm{y})=\sum_{\bm{k}\in\mathcal{K}_{K,L}}\varphi_{\bm{k}}(\bm{y}) \sum_{\bm{m}\in\mathbb{Z}^n\setminus\{\bm{0}_n\}}\hat{{U}}_{(\bm{k}_{\mathrm{I}}+2K\bm{m}_{\mathrm{I}},\bm{k}_{\mathrm{II}}+2L\bm{m}_{\mathrm{II}})},$$
where $\bm{0}_n$ denotes the $n$-dimensional zero vector. Let
$$g(\bm{k},\bm{m})=\|\bm{k}_{\mathrm{I}}+2K\bm{m}_{\mathrm{I}} +\bm{Q}(\bm{k}_{\mathrm{II}}+2L\bm{m}_{\mathrm{II}})\|^{2\alpha}+ \|\bm{k}_{\mathrm{II}}+2L\bm{m}_{\mathrm{II}}\|^{2\beta}.$$
Note that $g(\bm{k},\bm{m})=0$ if and only if $\bm{k}=\bm{m}=0$, and this is because the column vectors of the projection matrix $\bm{P}$ are $\mathbb{Q}$-linearly independent.
Using the Cauchy-Schwarz inequality, we have
\begin{equation}\label{eq:P-I}
\begin{aligned}
&\left\|\mathcal{P}_{K,L}u-\mathcal{I}_{K,L}u\right\|_{\mathcal{L}^2}^2\\
=&\sum_{\bm{k}\in\mathcal{K}_{K,L}} \bigg|\sum_{\bm{m}\in\mathbb{Z}^n\setminus\{\bm{0}_n\}}\hat{{U}}_{(\bm{k}_{\mathrm{I}}+2K\bm{m}_{\mathrm{I}},\bm{k}_{\mathrm{II}}+2L\bm{m}_{\mathrm{II}})}\bigg|^2\\
=&\sum_{\bm{k}\in\mathcal{K}_{K,L}}\bigg|\sum_{\bm{m}\in\mathbb{Z}^n\setminus\{\bm{0}_n\}}g(\bm{k},\bm{m})^{-1/2} g(\bm{k},\bm{m})^{1/2} \hat{{U}}_{(\bm{k}_{\mathrm{I}}+2K\bm{m}_{\mathrm{I}},\bm{k}_{\mathrm{II}}+2L\bm{m}_{\mathrm{II}})}\bigg|^2\\
\leq&\sum_{\bm{k}\in\mathcal{K}_{K,L}}\sum_{\bm{m}\in\mathbb{Z}^n\setminus\{\bm{0}_n\}} g(\bm{k},\bm{m})^{-1} \sum_{\bm{l}\in\mathbb{Z}^n\setminus\{\bm{0}_n\}}g(\bm{k},\bm{l}) |\hat{{U}}_{(\bm{k}_{\mathrm{I}}+2K\bm{l}_{\mathrm{I}},\bm{k}_{\mathrm{II}}+2L\bm{l}_{\mathrm{II}})}|^2.
\end{aligned}
\end{equation}

Next, we consider the series
$$\begin{aligned}
\sum_{\bm{m}\in\mathbb{Z}^n\setminus\{\bm{0}_n\}} g(\bm{k},\bm{m})^{-1}
=\sum_{\substack{\bm{m}_{\mathrm{I}}\in\mathbb{Z}_d\\ \bm{m}_{\mathrm{II}}\neq\bm{0}_{n-d}}} g(\bm{k},\bm{m})^{-1}+\sum_{\substack{\bm{m}_{\mathrm{I}}\neq\bm{0}_d\\\bm{m}_{\mathrm{II}}=\bm{0}_{n-d} }} g(\bm{k},\bm{m})^{-1}=b_1+b_2.
\end{aligned}$$
Since $\bm{k}_{\mathrm{II}}\in[-L,1-L,\cdots,L)^{n-d}$, we have
$$\|\bm{k}_{\mathrm{II}}+2L\bm{m}_{\mathrm{II}}\|^{2\beta}\geq L^{2\beta}\|\bm{m}_{\mathrm{II}}\|^{2\beta}.$$
For fixed $\bm{m}_{\mathrm{II}}$, there exists a unique $\bm{m}_{\mathrm{I}}'$ such that
$$\bm{k}_{\mathrm{I}}+2K\bm{m}_{\mathrm{I}}'+\bm{Q}(\bm{k}_{\mathrm{II}}+2L\bm{m}_{\mathrm{II}}) \in[-K,1-K,\cdots,K)^d,$$
and it is obvious $\bm{m}_{\mathrm{I}}'=\bm{0}_d$ if $\bm{m}_{\mathrm{II}}=\bm{0}_{n-d}$. Then,
$$\begin{aligned}
&\|\bm{k}_{\mathrm{I}}+2K\bm{m}_{\mathrm{I}}+ \bm{Q}(\bm{k}_{\mathrm{II}}+2L\bm{m}_{\mathrm{II}})\|^{2\alpha}\\
=&\|2K(\bm{m}_{\mathrm{I}}-\bm{m}_{\mathrm{I}}')+\bm{k}_{\mathrm{I}}+2K\bm{m}_{\mathrm{I}}'+ \bm{Q}(\bm{k}_{\mathrm{II}}+2L\bm{m}_{\mathrm{II}})\|^{2\alpha}\\
\geq&K^{2\alpha}\|(\bm{m}_{\mathrm{I}}-\bm{m}_{\mathrm{I}}')\|^{2\alpha}.
\end{aligned}$$
Thus,
$$\begin{aligned}
b_1=&\sum_{\bm{m}_{\mathrm{I}}\in\mathbb{Z}^d}~\sum_{\bm{m}_{\mathrm{II}}\neq\bm{0}_{n-d}} \big(\|\bm{k}_{\mathrm{I}}+2K\bm{m}_{\mathrm{I}} +\bm{Q}(\bm{k}_{\mathrm{II}}+2L\bm{m}_{\mathrm{II}})\|^{2\alpha}+ \|\bm{k}_{\mathrm{II}}+2L\bm{m}_{\mathrm{II}}\|^{2\beta}\big)^{-1}\\
\leq&\sum_{\bm{m}_{\mathrm{I}}\in\mathbb{Z}^d}~\sum_{\bm{m}_{\mathrm{II}}\neq\bm{0}_{n-d}} \big(K^{2\alpha}\|(\bm{m}_{\mathrm{I}}-\bm{m}_{\mathrm{I}}')\|^{2\alpha}+L^{2\beta}\|\bm{m}_{\mathrm{II}}\|^{2\beta}\big)^{-1}\\
\leq&~\max\big\{K^{-2\alpha},L^{-2\beta}\big\}\sum_{\bm{m}_{\mathrm{I}}\in\mathbb{Z}^d}~\sum_{\bm{m}_{\mathrm{II}}\neq\bm{0}_{n-d}} \big(\|\bm{m}_{\mathrm{I}}\|^{2\alpha}+\|\bm{m}_{\mathrm{II}}\|^{2\beta}\big)^{-1},\\
b_2\leq& \sum_{ \bm{m}_{\mathrm{I}}\neq\bm{0}_d}\|\bm{k}_{\mathrm{I}}+2K\bm{m}_{\mathrm{I}}+ \bm{Qk}_{\mathrm{II}}\|^{-2\alpha} 
\leq K^{-2\alpha}\sum_{\bm{m}_{\mathrm{I}}\neq\bm{0}_d}\|\bm{m}_{\mathrm{I}}\|^{-2\alpha}.
\end{aligned}$$
Combining the above two inequalities, we have
\begin{equation*}
    \sum_{\bm{m}\in\mathbb{Z}^n\setminus\{\bm{0}_n\}} g(\bm{k},\bm{m})^{-1}\lesssim\big(K^{-2\alpha}+L^{-2\beta}\big)
    \sum_{\bm{m}_{\mathrm{I}}\neq\bm{0}_d}~\sum_{\bm{m}_{\mathrm{II}}\neq\bm{0}_{n-d}} \big(\|\bm{m}_{\mathrm{I}}\|^{2\alpha}+\|\bm{m}_{\mathrm{II}}\|^{2\beta}\big)^{-1}.
\end{equation*}
The sufficient condition for the convergence series 
$\sum\limits_{\bm{m}_{\mathrm{I}}\neq\bm{0}_d
}\|\bm{m}_{\mathrm{I}}\|^{-2\alpha}
\sum\limits_{\bm{m}_{\mathrm{II}}\neq\bm{0}_{n-d}}\|\bm{m}_{\mathrm{II}}\|^{-2\beta}$, $\alpha>d/2$ 
and $\beta>(n-d)/2$, can be easily obtained by Hölder's inequality. Moreover, the condition $d/2\alpha+(n-d)/2\beta<1$ implies that there exist $\alpha'>d$ and $\beta'>n-d$ such that
$$\frac{\alpha'}{2\alpha}+\frac{\beta'}{2\beta}=1,$$
and by Young's inequality, we have
$$\sum_{\bm{m}_{\mathrm{I}}\neq\bm{0}_d}~\sum_{\bm{m}_{\mathrm{II}}\neq\bm{0}_{n-d}} \big(\|\bm{m}_{\mathrm{I}}\|^{2\alpha}+\|\bm{m}_{\mathrm{II}}\|^{2\beta}\big)^{-1} \leq\sum_{\bm{m}_{\mathrm{I}}\neq\bm{0}_d}~\sum_{\bm{m}_{\mathrm{II}}\neq\bm{0}_{n-d}} \|\bm{m}_{\mathrm{I}}\|^{-\alpha'}\|\bm{m}_{\mathrm{II}}\|^{-\beta'},$$
which is obviously convergent. The convergence of the series yields
$$\sum_{\bm{m}\in\mathbb{Z}^n\setminus\{\bm{0}_n\}} g(\bm{k},\bm{m})^{-1} \lesssim\big(K^{-2\alpha}+L^{-2\beta}\big).$$
Therefore, it follows from \cref{eq:P-I} that
$$\begin{aligned}
\left\|\mathcal{P}_{K,L}u-\mathcal{I}_{K,L}u\right\|_{\mathcal{L}^2}^2
\lesssim &\big(K^{-2\alpha}+L^{-2\beta}\big)\sum_{\bm{k}\in\mathcal{K}_{K,L}} \sum_{\bm{l}\in\mathbb{Z}^n\setminus\{\bm{0}_n\}}g(\bm{k},\bm{l}) |\hat{{U}}_{(\bm{k}_{\mathrm{I}}+2K\bm{l}_{\mathrm{I}},\bm{k}_{\mathrm{II}}+2L\bm{l}_{\mathrm{II}})}|^2\\
\leq &\big(K^{-2\alpha}+L^{-2\beta}\big)|u|_{\alpha,\beta}^2.
\end{aligned}$$

Moreover,
$$\begin{aligned}
\left|\mathcal{P}_{K,L}u-\mathcal{I}_{K,L}u\right|_{\mu,\nu}^2=&\sum_{\bm{k}\in\mathcal{K}_{K,L}} \big(\|\bm{k}_{\mathrm{I}} +\bm{Qk}_{\mathrm{II}}\|^{2\mu}+ \|\bm{k}_{\mathrm{II}}\|^{2\nu}\big) |\hat{{U}}_{\bm{k}}-\tilde{u}_{\bm{k}}|^2\\
\leq&\big(K^{2\mu}+L^{2\nu}\big)\left\|\mathcal{P}_{K,L}u-\mathcal{I}_{K,L}u\right\|_{\mathcal{L}^2}^2\\
\lesssim&\big(K^{-2\alpha}+L^{-2\beta}\big)\big(K^{2\mu}+L^{2\nu}\big)|u|_{\alpha,\beta}^2.
\end{aligned}$$
Thus, from \Cref{cr:proj}, we have
$$\begin{aligned}
\left|u-\mathcal{I}_{K,L}u\right|_{\mu,\nu}^2\leq&\left|u-\mathcal{P}_{K,L}u\right|_{\mu,\nu}^2 +\left|\mathcal{P}_{K,L}u-\mathcal{I}_{K,L}u\right|_{\mu,\nu}^2\\
\lesssim&\big(K^{-2\alpha}+L^{-2\beta}\big)\big(K^{2\mu}+L^{2\nu}\big)|u|_{\alpha,\beta}^2.
\end{aligned}$$
Similarly, we can also obtain
$$\|u-\mathcal{I}_{K,L}u\|_{\mu,\nu}^2\lesssim\big(K^{-2\alpha}+L^{-2\beta}\big)\big(K^{2\mu}+L^{2\nu}\big)\|u\|_{\alpha,\beta}^2.$$
The original proposition is established.
\end{proof}

\end{document}

%% file: ex_shared.tex
\headers{IWFPM and application to QSEs}{K. Jiang, X. Li, Y. Ma, J. Zhang, P. Zhang, Q. Zhou}

\title{Irrational-window-filter projection method and application to quasiperiodic Schr\"odinger eigenproblems
\thanks{Submitted to xxx.
  \funding{This work was partially supported by the National Key R\&D Program of China (2023YFA1008802), the NSFC grants (12171412, 12288101), the Science and Technology Innovation Program of Hunan Province (2024RC1052).
  XL was supported partially by Scientific Research Fund of Hunan Provincial Education Department (23A0127). 
  We are also grateful to the High Performance Computing Platform of Xiangtan University.
 }
}}

\author{Kai Jiang\thanks{
		Hunan Key Laboratory for Computation and Simulation in Science and Engineering,
		Key Laboratory of Intelligent Computing and Information Processing of Ministry
		of Education, School of
		Mathematics and Computational Science, Xiangtan University, Xiangtan, Hunan,
		China, 411105.	
		(\email{kaijiang@xtu.edu.cn},
	\email{lixy1217@xtu.edu.cn}, \email{mayao@smail.xtu.edu.cn}, 		\email{zhangjuan@xtu.edu.cn},  \email{qizhou@smail.xtu.edu.cn}.)
}
	\and Xueyang Li\footnotemark[2] 
	\and Yao Ma\footnotemark[2]
	\and Juan Zhang\footnotemark[2]
	\and Pingwen  Zhang\thanks{
		School of Mathematics and Statistics, Wuhan University, Wuhan, China, 430072; 
		School of Mathematical Sciences,
		Peking University, Beijing, China, 100871.
	(\email{pzhang@pku.edu.cn})	}
	\and
	Qi Zhou\footnotemark[2]
	}

\usepackage{mathrsfs,amsmath,amssymb,bm}
\usepackage{lipsum}
\usepackage{amsfonts}
\usepackage{caption,graphicx, subfigure}
\usepackage{epstopdf}
\usepackage{makecell, rotating, bbding}
\usepackage{multirow}
\usepackage{algorithmic, algorithm}
\usepackage{enumerate}
\usepackage{booktabs}
\usepackage[misc]{ifsym}
\usepackage{mathtools}
\usepackage{amsmath}

\newtheorem{thm}{Theorem}[section]

\newtheorem{remark}[thm]{Remark}
\newtheorem{example}[thm]{Example}

\DeclareMathOperator*{\argmin}{\mathrm{argmin}}

\usepackage[draft]{changes}
\usepackage{lipsum}

\newsiamremark{hypothesis}{Hypothesis}
\crefname{hypothesis}{Hypothesis}{Hypotheses}
\newsiamthm{claim}{Claim}

\usepackage{amsopn}
\DeclareMathOperator{\diag}{diag}

\makeatletter
\newcommand*{\addFileDependency}[1]{
  \typeout{(#1)}
  \@addtofilelist{#1}
  \IfFileExists{#1}{}{\typeout{No file #1.}}
}
\makeatother
